\documentclass{IEEEtran}
\usepackage{cite}
\usepackage{amsmath,amssymb,amsfonts}

\usepackage{graphicx}
\usepackage{textcomp}
\def\BibTeX{{\rm B\kern-.05em{\sc i\kern-.025em b}\kern-.08em
    T\kern-.1667em\lower.7ex\hbox{E}\kern-.125emX}}

\newcommand{\hide}[1]{}
\usepackage{physics}
\newcommand*{\QEDB}{\hfill\ensuremath{\blacksquare}}%
\usepackage{enumitem}

\usepackage{algorithm} 
\usepackage{algpseudocode} 

\usepackage[english]{babel}
\usepackage{amsthm}
\usepackage{xcolor}

\theoremstyle{remark}

\newtheorem{theorem}{Theorem}
\newtheorem{lemma}{Lemma}
\newtheorem{remark}{Remark}
\newtheorem{assumption}{Assumption}
\newtheorem{corollary}{Corollary}

\def\mz{\mu_z}
\def\mt{\mu_\theta}

\def \sigminwJ{\sigma^{\overline{J}}_{\mathrm{min}}}
\def \sigmaxLs{\sigma_{\mathrm{max}}^{L_s}}

\def\va{v_\alpha}
\def\vas{v_\alpha^\star}

\usepackage{tikz}
\usetikzlibrary{decorations.pathreplacing,calc}
\newcommand{\tikzmark}[1]{\tikz[overlay,remember picture] \node (#1) {};}

\newcommand*{\AddNote}[4]{%
    \begin{tikzpicture}[overlay, remember picture]
        \draw [decoration={brace,amplitude=0.5em},decorate,ultra thick,black]
            ($(#3)!(#1.north)!($(#3)-(0,1)$)$) --  
            ($(#3)!(#2.south)!($(#3)-(0,1)$)$)
                node [align=center, text width=1.5cm, pos=0.5, anchor=west] {#4};
    \end{tikzpicture}
}%

\usepackage{multicol}

\usepackage{subcaption}
\usepackage{caption}

\usepackage{mathtools}

\usepackage{booktabs}
\usepackage{siunitx}

\newcommand{\ceil}[1]{\left\lceil {#1} \right\rceil}

\begin{document}

\title{\huge Communication Efficient Curvature Aided Primal-dual
Algorithms for Decentralized Optimization
}


\author{Yichuan Li, Petros G. Voulgaris, Du\v{s}an M. Stipanovi\'{c}, and Nikolaos M. Freris

\thanks{Yichuan Li is with the Coordinated Science Laboratory and the Department of Mechanical Science and Engineering, University of Illinois Urbana-Champaign, IL 61820 USA (email: yli129@illinois.edu).}

\thanks{Petros G. Voulgaris is with the Department of Mechanical Engineering, University of Nevada, Reno, NV 89557, USA (email: pvoulgaris@unr.edu).}

\thanks{Du\v{s}an M. Stipanovi\'{c} is with the Coordinated Science Laboratory and the Department of Industrial and Enterprise Systems Engineering, University of Illinois Urbana-Champaign, IL 61820 USA (email: dusan@illinois.edu).}

\thanks{Nikolaos M. Freris is with the School of Computer Science, University of Science and Technology of China, Hefei, Anhui, 230027, China (email: nfr@ustc.edu.cn).} 

\thanks{Freris (correspondence) was supported by the Ministry of Science and Technology of China under grant 2019YFB2102200. Voulgaris and Li were supported by NSF grants CCF-1717154, CPS-1932529 and CMMI-2137764. }
}

\maketitle

\begin{abstract}
This paper presents a family of algorithms for decentralized convex composite problems. We consider the setting of a network of agents that cooperatively minimize a global objective function composed of a sum of local functions plus a regularizer. Through the use of intermediate consensus variables, we remove the need for inner communication loops between agents when computing curvature-guided updates. A general scheme is presented which unifies the analysis for a plethora of computing choices, including gradient descent, Newton updates, and BFGS updates. Our analysis establishes sublinear convergence rates under convex objective functions with Lipschitz continuous gradients, as well as linear convergence rates when the local functions are further assumed to be strongly convex. Moreover, we explicitly characterize the acceleration due to curvature information. Last but not the least, we present an asynchronous implementation for the proposed algorithms, which removes the need for a central clock, with linear convergence rates established in expectation under strongly convex objectives. We ascertain the effectiveness of the proposed methods with numerical experiments on benchmark datasets. 
\end{abstract}


\section{Introduction}
\label{sec1}
\IEEEPARstart{T}{he} proliferation of mobile devices with computation and communication capabilities has fueled the surge of applications of distributed optimization in various fields. Examples include distributed control\hide{\cite{discontrol2009,discontrol2013}}, wireless sensor networks\hide{\cite{sensor2004}}\hide{\cite{sensor2013}}, power grid management\hide{\cite{power2013,power2016}}, and large-scale machine learning\hide{\cite{ML2012}\cite{ML2018}}\cite{discontrol2018,sensor2011,power2020,ML2020,new1}. A canonical problem in distributed optimization assumes a network of agents collaboratively optimizing a global objective function through message passing with immediate neighbors. In specific, we consider the following optimization problem:
\begin{gather}
    \underset{\hat{x}\in \mathbb{R}^d}{\text{minimize}}\quad \left\{\sum_{i=1}^m f_i(\hat{x})+g(\hat{x})\right\}, \label{prob1}
\end{gather}
where each $f_i(\cdot):\mathbb{R}^d \to \mathbb{R}$ is a convex and smooth function accessible only by agent $i$ while $g(\cdot):\mathbb{R}^d\to \mathbb{R}$ is a convex (possibly nonsmooth) regularizer. The inclusion of the regularizer is multi-faceted, e.g., it serves for promoting desired structures in the decision vector, such as sparsity in controller designs using the $\ell_1$-norm \cite{controllerl2013}, preventing overfitting in machine learning using the squared $\ell_2$-norm \cite{ridge1998}, and enforcing constraints using indicator functions of convex sets. 

First-order methods \cite{1storder2009,1st2015,pgextra,1st2018,p2d2,1st2019,1st2021,new2} using (sub)-gradient information constitute popular choices for solving (\ref{prob1}) due to their economical computational costs and simple implementation. In first-order methods, agents compute updates by using local gradients combined with averaged information from their neighbors. For the case with no regularizer ($g(\cdot)\equiv 0$), \cite{1st2015,1st2018,1st2019} exploit the history of gradient and iterate values to achieve linear convergence rate for strongly convex objectives. \cite{1st2021} provides a unified framework for designing various first-order schemes for the general problem (\ref{prob1}). When nonsmooth regularizers are present, existing work almost exclusively applies proximal gradient type of updates: each agent first performs gradient descent on the smooth part of the objective function and then invokes the proximal operator associated with the nonsmooth regularizer $g(\cdot)$. Nonetheless, using only first-order information suffers from slow convergence speed and thus requires a large number of total iterations to reach a prescribed accuracy. This constitutes a key limitation for first-order methods, which is most pronounced in applications where high-accuracy solutions are pursued in a few rounds of iterations, for example, due to high communication costs.

A natural option for accelerating the convergence is to use second-order information for local updates. Most second-order methods \cite{2nd2013,2nd2014,2nd2017,2nd2020} for solving (\ref{prob1}) focus on cases where the objective function is smooth, i.e, $g(\cdot)=0$. One reason is that even when proximal gradient steps are efficiently computable, proximal Newton steps require significantly more computational resources due to the Hessian scaling in the evaluation of the proximal operator. Another challenge in designing second-order methods lies in constructing distributedly computable Newton updates. Computing curvature-guided updates requires solving a linear system that, in general, involves global information, whence a direct application of the Newton method is not feasible. Moreover, the standard Newton method requires backtracking line search to select appropriate step sizes for ensuring global convergence \cite{Numerical}. Such operations incur heavy communication burdens in the form of collecting all local objective function values in the network; this necessitates extensive message passing between agents or the presence of a centralized coordinator. Authors in \cite{2nd2014} propose to use matrix splitting techniques in the dual problem, so that the Hessian inverse admits a distributedly computable Taylor expansion. By truncating the Taylor series to $K$ terms, agents may compute local updates with an additional $K$ rounds of communication loops with their neighbors. With $g(\cdot)=0$, \cite{2nd2017} and \cite{2nd2020} use similar matrix splitting techniques to solve a penalized version of (\ref{prob1}) where the former presents a synchronous scheme and the latter extends it to asynchronous settings. We note that \cite{2nd2017} and \cite{2nd2020} are effectively solving a different problem (penalized version) compared to (\ref{prob1}) when using constant stepsizes, and therefore do not converge to the exact solution.

Another popular line of algorithmic design for solving (\ref{prob1}) is based on primal-dual methods, such as the Generalized Method of Multipliers, the Augmented Lagrangian Method, and the Alternating Direction Method of Multipliers (ADMM) \cite{Bert1989,Boyd2011}. In the setting of distributed primal-dual algorithms \cite{Mateos2010, Duchi2012,ermin2012, Wei2014, Hong2017,Latafat2019}, agents solve a sub-optimization problem at each iteration, which often involves multiple inner loops and thus induces heavy computation burden. Several approximation schemes \cite{Jakovetic2015,Ling2015,chang2015,DQM2016,Mokhtari2016,Eisen2019} were proposed to replace the exact minimization step with one or multiple update steps using approximated models of the augmented Lagrangian. It has been shown that by appropriately choosing the mixing matrices and the augmented Lagrangian model, primal-dual algorithms can recover several accelerated primal-only algorithms using gradient and iterate tracking techniques\cite{jakovetic2020}. Further acceleration can be achieved by resorting to Newton or quasi-Newton primal updates \cite{DQM2016,Mokhtari2016,Eisen2019}. However, all of them are synchronous algorithms considering the smooth problem ($g(\cdot)\equiv 0$) and  \cite{Mokhtari2016}, \cite{Eisen2019}  require multiple inner communication loops at each iteration of the algorithm. In such scenarios, despite improving the convergence speed, it is not clear whether the overall communication costs can be reduced due to the additional communication rounds per iteration. In emerging applications such as multi-agent Cyberphysical Systems \cite{CPS2012} and Federated Learning \cite{FL2017,fedadmm}, high responsiveness and reducing communication costs are of primordial importance. This motivates the development of methods with accelerated convergence as well as with guaranteed low communication costs, which is the focus of this paper.

\noindent\textbf{Contributions}:
\begin{itemize}[noitemsep,topsep=0pt]
    \item We introduce a framework for designing distributed primal-dual algorithms for (\ref{prob1}) with a nonsmooth regularization function. Through the use of intermediate consensus variables, we decouple the primal subproblem pertaining to an agent from those of its neighbors. As a result, we obtain a block-diagonal Hessian that allows us to incorporate curvature information in local updates \emph{without additional communication}. This is in contradistinction with the state-of-the-art, where multiple communication inner loops are required to compute (quasi) Newton updates.
    \item Using this framework, we propose DistRibuted cUrvature aided prImal Dual algorithms (DRUID), a family of algorithms that offer flexible choices of updating schemes, including gradient, Newton, and BFGS type of updates. Furthermore, we present a unified analysis framework for this class of algorithms, which not only establishes $\mathcal{O}(\frac{1}{T})$ convergence rate to optimality 
    under convex objectives, but also theoretically reveals the discrepancies among them. When strong convexity is further assumed, we establish linear convergence rates for this class of algorithms, and once again quantify the acceleration.
    \item We devise an asynchronous extension for this class of algorithms, and establish linear convergence rates in expectation, under strong convexity. This setting removes the need for a central clock in the network, and further allows for an arbitrary number of agents to be active at each iteration. We demonstrate the merits of the proposed framework through simulations using real-life datasets. 
\end{itemize}
\noindent \textbf{Notation}: We represent column vectors $x\in \mathbb{R}^d$ with lower case letters, matrices $A\in \mathbb{R}^{n\times m}$ with capital letters, and matrix transpose as $A^\top$. We also use $[A,B]$ and $[A;B]$ to respectively denote row and column stacking (for matrices with equal numbers of rows or columns, respectively). Superscript denotes the sequence index while subscript denotes the vector component. For example, $x_i^t$ represents the vector component held by agent $i$ at iteration $t$. Moreover, $[A]_{ij}$ denotes the $ij$-th entry of matrix $A$. If a norm specification is not provided, $\norm{x}$ and $\norm{A}$ represent the vector Euclidean norm and the induced matrix norm, respectively. For a positive definite matrix $P\succ 0$, we define $\norm{x}_P:=\sqrt{x^\top P x}$. The set $\{1,\dots,m\}$ is abbreviated as $[m]$ and the proximal mapping associated with a function $g(\cdot):\mathbb{R}^d\to \mathbb{R}$ is defined as $\textbf{prox}_{\tfrac{g}{\mu}}(v) := \underset{\theta\in \mathbb{R}^d}{\text{argmin}}\left\{g(\theta)+\frac{\mu}{2}\norm{\theta-v}^2\right\}$. We further denote the identity matrix of dimension $d$ as $I_d$ and the Kronecker product between two matrices of arbitrary dimension $A,B$ as $A\otimes B$. 

\section{Preliminaries} 
In this section, we begin with reformulating problem (\ref{prob1}) to the consensus setting that is used for our development in Section \ref{reformulation}. \hide{We provide some background on using ADMM to solve the reformulated problem in Section \ref{background} and an introduction to quasi-Newton methods in Section \ref{intro_quasi}.} 
  
\subsection{Problem formulation} 
\label{reformulation}
We capture the network topology by an undirected graph $G=\{\mathcal{V},\mathcal{E}\}$ where $\mathcal{V}:=[m]$ denotes the vertex set and the edge set $\mathcal{E}\subseteq \mathcal{V}\times \mathcal{V}$ contains the pair $(i,j)$ if and only if agent $i$ can communicate with agent $j$. We do not consider self loops, i.e., $(i,i)\notin \mathcal{E}$ for any $i\in[m]$. For notational convenience, we enumerate the edge set (arbitrary order) and use $\mathcal{E}_k$ to denote the $k$-th edge, $k\in[n]$, where $n\coloneqq \abs{\mathcal{E}}$ is the number of edges. Moreover, the set of neighbors of agent $i$ is defined as $\mathcal{N}_i:=\{j\in \mathcal{V}:(i,j)\in \mathcal{E}\}$. Using the above definitions, we reformulate problem (\ref{prob1}) to the following consensus formulation, by introducing local decision variables $x_i$ at corresponding agent $i$, as well as edge variables $z_{ij}$ for $(i,j)\in \mathcal{E}$. The consensus formulation is given by: 
\begin{equation}
\begin{aligned}
    \underset{x_i,\theta,z_{ij}\in\mathbb{R}^d}{\mathrm{minimize}}&\quad \left\{\sum_{i=1}^m f_i(x_i)+g(\theta)\right\},\\
    \text{s.t.}\,\,x_i=z_{ij}&= x_j,\,\forall\,i\in [m]\text{ and }j\in \mathcal{N}_i.\\
    x_l &=\theta,\,\text{for one arbitrary $l\in [m]$}.
\end{aligned}\label{prob2}
\end{equation}
Note that we have also introduced $\theta$ to separate the argument of the smooth and nonsmooth functions and only enforce the equality constraint for $\theta$ at the $l$-th agent as $x_l=\theta$, where $l$ can be arbitrarily selected. We emphasize that this agent is not a central coordinator, but rather the agent whose local updates factor in the nonsmooth regularizer. This is without loss of generality and induces minimal computational overhead from evaluating proximal mappings. Assuming $G$ is connected, it is easy to check that (\ref{prob2}) is equivalent to (\ref{prob1}) since their optima coincide, i.e., $\hat{x}^\star= x_i^\star=z_{ij}^\star=\theta^\star$, $\forall\,i\in[m]$ and $j\in \mathcal{N}_i$. This is achieved by satisfying the consensus constraints in (\ref{prob2}). We note that consensus can be enforced by simply letting $x_i=x_j$, i.e., without intermediate consensus variables $\{z_{ij}\}$. However, the introduction of intermediate variables is key to our design: the purpose of $\{z_{ij}\}$ is to decouple $x_i$ from its neighbors so that we achieve a block-diagonal Hessian for the augmented Lagrangian. A block-diagonal Hessian allows agents to compute the (quasi) Newton steps \emph{without additional communication with their neighbors}. We provide further discussion on this choice in Section \ref{section3}. 

We proceed to define the source and destination matrices $\hat{A}_s,\hat{A}_d\in \mathbb{R}^{n\times m}$. Each row of $\hat{A}_s$ and $\hat{A}_d$ corresponds to an edge $\mathcal{E}_k$ in the graph, $k\in [n]$: $[\hat{A}_s]_{ki}=[\hat{A}_d]_{kj}=1$ if and only if $\mathcal{E}_k=(i,j)$, and $0$ otherwise. Problem (\ref{prob2}) can then be compactly expressed using the concatenated column vectors $x:=[x_1^\top,\dots,x_m^\top]^\top,z:=[z_1^\top,\dots,z_n^\top]^\top$ (we note a slight abuse of notation in using $z_k\equiv z_{ij}$ where $k\in[n]$ is the corresponding edge $(i,j)\in\mathcal{E}$ in the enumeration order) as:   
\begin{equation}
\begin{aligned}
    \underset{x\in \mathbb{R}^{md},\theta \in \mathbb{R}^d,z\in \mathbb{R}^{nd}}{\text{minimize}}\quad  & \Big \{F(x)+g(\theta)\Big\},\\
    \text{s.t.}\,\, Ax = \begin{bmatrix}
    \hat{A}_s \otimes I_d \\ 
    \hat{A}_d \otimes I_d\end{bmatrix}x &= 
    \begin{bmatrix}
    I_{nd} \\ I_{nd}
    \end{bmatrix}z= Bz,  \\ 
    S^\top x&= \theta ,
\end{aligned}\label{prob3}
\end{equation}
where $F(x):= \sum_{i=1}^m f_i(x_i)$ and matrices $A$ and $B$ are obtained by stacking the matrices as shown in (\ref{prob3}). We further define $S:= (s_l \otimes I_d)\in \mathbb{R}^{md\times d}$ where $s_l\in \mathbb{R}^m$ is an all-zero vector except for the $l$-th entry being one. In other words, the $S^\top$ matrix serves to select the $l$-th component of $x$ held by the agent $l$, i.e., $S^\top x = x_l$. We proceed to present some identities that associate source and destination matrices to the incidence and Laplacian matrices corresponding to the graph topology in the following. 
\begin{subequations}
\begin{align}
    \hat{E}_s &= \hat{A}_s -\hat{A}_d,\,\, \hat{E}_u=\hat{A}_s+\hat{A}_d, \label{4a}\\
    \hat{L}_s &= \hat{E}_s^\top \hat{E}_s,\,\,\hat{L}_u= \hat{E}_u^\top \hat{E}_u,\label{4b}\\
    \hat{D} &= \tfrac{1}{2}(\hat{L}_s+\hat{L}_u)= \hat{A}_s^\top \hat{A}_s+\hat{A}_d^\top \hat{A}_d, \label{4c}
\end{align}
\end{subequations}
where $\hat{E}_s,\hat{E}_u\in \mathbb{R}^{n\times m}$ are signed and unsigned graph incidence matrices and $\hat{L}_s,\hat{L}_u\in \mathbb{R}^{m\times m}$ are signed and unsigned graph Laplacian matrices respectively. The diagonal matrix $\hat{D}\in\mathbb{R}^{m\times m}$ denotes the graph degree matrix with entries $D_{ii}=\abs{\mathcal{N}_i}$. We further introduce the block extensions to the dimension $d$, that is $E_s:=\hat{E}_s\otimes I_d$ and similarly for $E_u,L_s,L_u$, and $D$. \hide{, i.e., the Kronecker product between the matrix and $I_d$ is denoted without the hat.}  
\subsection{Background on ADMM} \label{background}
We begin by defining the augmented Lagrangian for problem (\ref{prob3}):
\begin{gather}
    \mathcal{L}(x,\theta,z;y,\lambda) := F(x)+g(\theta)+y^\top (Ax-Bz)\nonumber\\
    +\lambda^\top (S^\top x-\theta)
    +\tfrac{\mu_z}{2}\norm{Ax-Bz}^2+\tfrac{\mu_\theta}{2}\norm{S^\top x-\theta}^2, \label{AL}
\end{gather}
where $y\in\mathbb{R}^{2nd},\lambda\in\mathbb{R}^{d}$ are Lagrange multipliers associated with the constraints $Ax=Bz$ and $S^\top x=\theta$, respectively. Note that since penalty coefficients of the quadratic terms are closely related to dual step sizes, we have separated them into $\mz$ and $\mt$ to offer broader choices of selection. ADMM solves (\ref{prob3}), equivalently (\ref{prob2}) and (\ref{prob1}), by sequentially minimizing the augmented Lagrangian (\ref{AL}) over each of the primal variables $(x,\theta,z)$, and then performs gradient ascent on the dual variables $(y,\lambda)$:
\begin{subequations}
\begin{align}
    x^{t+1}&= \underset{x}{\text{argmin}}\,\,\mathcal{L} (x,\theta^t,z^t;y^t,\lambda^t) ,\label{admm1}\\
    \theta^{t+1} &= \underset{\theta}{\text{argmin}}\,\,\mathcal{L} (x^{t+1},\theta,z^t;y^t,\lambda^t) ,\label{admm2}\\
    z^{t+1} & = \underset{z}{\text{argmin}}\,\,\mathcal{L} (x^{t+1},\theta^{t+1},z;y^t,\lambda^t) ,\label{admm3}\\
    y^{t+1} &=y^t+\mu_z (Ax^{t+1}-Bz^{t+1}),\label{admm4}\\
    \lambda^{t+1} &=\lambda^t+\mu_\theta(S^\top x^{t+1}-\theta^{t+1})\label{admm5}.
\end{align}
\end{subequations}
The above iterations fall into the category of 3-block ADMM which is not guaranteed to converge for arbitrary $\mu_z,\mu_\theta>0$ \cite{3block}. Step (\ref{admm1}) requires a solution to a sub-optimization problem which often involves multiple inner-loops for general objective functions, and therefore becomes the most expensive step in ADMM. Executing step (\ref{admm2}) bears the complexity of computing the proximal mapping of  the regularization function $g(\cdot)$. For commonly used $g(\cdot)$, such as the $\ell_1$-norm, squared $\ell_2$-norm, and indicators of several convex sets, a closed-form solution exists. For other cases, one would often resort to the fact that the proximal operator is separable, Lipschitz continuous with constant 1, firmly nonexpansive, and the associated Moreau envelope function is continuously differentiable irrespective of the function $g(\cdot)$, to devise efficient algorithms to approximate the proximal mapping. We refer readers to \cite{proximal2014} for more details. Step (\ref{admm3}) results from the introduction of \hide{consensus variables}$\{z_{ij}\}$-variables, and it does not require explicit evaluation as we demonstrate in the Section \ref{section3}. 

\subsection{Introduction to quasi-Newton methods}\label{intro_quasi}
Quasi-Newton methods \cite{Numerical} constitute a class of methods that aim to accelerate convergence using curvature information of the objective function without solving a linear system as in the Newton method. Specifically, the update direction $u^t\in\mathbb{R}^d$ in quasi-Newton methods is given by:
\begin{gather*}
    u^t = (H^t)^{-1} \nabla f(x^t),
\end{gather*}
where $(H^t)^{-1} \succ 0$ is some matrix (the inverse is just notation for ease of exposition, and no inversion is needed) that approximates the Hessian inverse. One of the main advantages of quasi-Newton methods lies in the fact that $(H^t)^{-1}$ is explicitly available so computing $u^t$ amounts to performing matrix multiplication at the cost of $\mathcal{O}(d^2)$ for general problems, as compared to solving a linear system with computational costs $\mathcal{O}(d^3)$ in Newton method. Many schemes exist for estimating $(H^t)^{-1}$ and in subsequent discussions, we focus on the one proposed by Broyden, Fletcher, Goldfarb, and Shanno (BFGS) \cite{Broyden1970,Fletcher1970,Goldfarb1970,Shanno1970}, which is considered to be the most effective in terms of acceleration and self-correcting capabilities\cite{Numerical}.

We define the consecutive iterate and gradient differences as:
\begin{gather}
    s^t =x^{t+1}-x^t,\,\,\text{and}\,\, q^t =\nabla f(x^{t+1})-\nabla f(x^t).
\end{gather}
BFGS requires the updated Hessian inverse approximation $(H^{t+1})^{-1}$ to satisfy the following secant condition:
\begin{gather}
    (H^{t+1})^{-1}q^t = s^t,\label{secant}
\end{gather}
which is motivated by the fact that the exact Hessian inverse satisfies (\ref{secant}) as $x^{t+1}$ tends to $x^t$.\hide{To ensure positive definiteness of $(H^{t+1})^{-1}$, it must hold that $(q^t)^\top s^t > 0$, as can be seen by premultiplying (\ref{secant}) with $(q^t)^\top$. This is satisfied automatically with convex $f(\cdot)$.} However, the secant condition alone is not enough to specify $(H^{t+1})^{-1}$. BFGS proposes to select $(H^{t+1})^{-1}$ by further requiring the updated estimate to be close to the previous one in the following sense:
\begin{gather}
    \underset{H^{-1}}{\text{minimize}}\quad  \norm{H^{-1}-(H^t)^{-1}}_{W} \label{bfgs_prob}\\
    \text{s.t.} \,\,H^{-1} = (H^{-1})^\top,\quad H^{-1}q^t=s^t, \nonumber
\end{gather}
where $\norm{M}_{W}:=\norm{W^{\tfrac{1}{2}}MW^{\tfrac{1}{2}}}_F$ denotes the weighted Frobenius norm with $W$ being the average Hessian \cite{Numerical}. Problem (\ref{bfgs_prob}) admits a closed-form solution which gives rise to the update formula for $(H^{t+1})^{-1}$ in BFGS: 
\begin{gather}
    (H^{t+1})^{-1} = \left(I-\rho^ts^t(q^t)^\top\right)(H^t)^{-1}\left(I-\rho^t q^t(s^t)^\top\right)
    \nonumber\\
    +\rho^t s^t (s^t)^\top, \label{bfgs_formula}
\end{gather}
where $\rho^t= 1/(q^t)^\top s^t$. \hide{We note that an update formula for the estimated Hessian can also be derived by using the Sherman-Morrison formula as follows: 
\begin{gather}
    H^{t+1} = H^t - \frac{H^ts^t (s^t)^\top H^t}{(s^t)^\top H^t s^t}+\frac{q^t (q^t)^\top }{(q^t)^\top s^t}.\label{bfgs_formula_2} 
\end{gather}}
By using only gradient information as in first-order methods, BFGS iteratively constructs a Hessian estimate of the objective function as in (\ref{bfgs_formula}) that is accurate enough to achieve superlinear convergence rates. However, the direct application of BFGS does not admit a distributed implementation as can be seen in the formula (\ref{bfgs_formula}) where computing $s^t (q^t)^\top$ involves global operations and message passing among agents. In the following section, we introduce BFGS updates in the framework of primal-dual algorithms that are not only distributedly computable, but also \emph{retain the same communication costs} as their first-order counterparts.  

\section{Algorithmic development}
\label{section3}
An approximated augmented Lagrangian $\widehat{\mathcal{L}}(\cdot)$ can be obtained using second-order expansion as follows: $$ 
    \widehat{\mathcal{L}} (x,\theta^t,z^t;y^t,\lambda^t)= \mathcal{L}^t(x^t)+ (x-x^t)^\top\nabla_x\mathcal{L}^t
    +\tfrac{1}{2}\norm{x-x^t}_{H^t}^2, $$
where we abbreviated the $\mathcal{L}(x^t,\theta^t,z^t;y^t,\lambda^t)$ as $\mathcal{L}^t(x^t)$, and the selection of $H^t$ is a means for designing a range of methods as will be subsequently elaborated in this section. We obtain a closed form solution when minimizing $\widehat{\mathcal{L}}(\cdot)$ over $x$ and replace step (\ref{admm1}) with the following one-step update:
\begin{gather}
    x^{t+1}= x^t- (H^t)^{-1}\nabla_x \mathcal{L}(x^t,\theta^t,z^t;y^t,\lambda^t) \label{primal}.
\end{gather}
By completion of squares, step (\ref{admm2}) admits an analytical expression through the proximal mapping:
\begin{gather}
    \theta^{t+1}= \textbf{prox}_{g/\mt}(S^\top x^{t+1}+\tfrac{1}{\mt}\lambda^t).\label{prox}
\end{gather}
Moreover, since the augmented Lagrangian is quadratic with respect to $z$, it follows that $z^{t+1}$ of step (\ref{admm3}) can be computed by solving the following linear system of equations:
\begin{gather}
    B^\top y^t+\mu_z B^\top(Ax^{t+1}-Bz^{t+1})= 0.\label{zupdate}
\end{gather}
Dual variables are updated in verbatim as in steps (\ref{admm4}) and (\ref{admm5}). We note that dual updates can be performed in parallel once primal updates are completed. Before we explicate the choice for $H^t$, we present a lemma that allows for efficient implementation of (\ref{primal})--(\ref{zupdate}) and (\ref{admm4})--(\ref{admm5}) under appropriate initialization.

\begin{lemma}\label{lemma1} Recall the identities in (\ref{4a})--(\ref{4c}) and the definitions thereafter. We express the dual variable as $y^t=[\alpha^t;\beta^t]$, $\alpha^t,\beta^t\in\mathbb{R}^{nd}$. If $y^0$ and $z^0$ are initialized so that
$
    \alpha^0+\beta^0=0$ and $
    z^0 = \tfrac{1}{2}E_u x^0
$, then $\alpha^t+\beta^t=0$ and $z^t=\tfrac{1}{2}E_ux^t$ for all $t\geq 0$. Moreover, defining $\phi^t=E_s^\top \alpha^t$, we equivalently express the updates (\ref{primal})--(\ref{zupdate}), (\ref{admm4})--(\ref{admm5}) as:
\begin{subequations}\label{updates}
\begin{align}
        x^{t+1}&=x^t - (H^t)^{-1}\big[\nabla F(x^t)+\phi^t+ S \lambda^t+\tfrac{\mz}{2}L_sx^t\nonumber  \\&+\mt S (S^\top x^t-\theta^t)\big],\label{lemma1_1}\\
        \theta^{t+1}&= \textbf{prox}_{g/\mu_\theta}(S^\top x^{t+1}+\tfrac{1}{\mu_\theta}\lambda^t),\label{lemma1_2}\\
       \phi^{t+1}&=\phi^t+\tfrac{\mu_z}{2}L_s x^{t+1} ,\label{lemma1_3}\\
        \lambda^{t+1}&=\lambda^t+\mu_\theta(S^\top x^{t+1}-\theta^{t+1}).\label{lemma1_4}
\end{align}
\end{subequations}
\textit{Proof}: See Appendix \ref{AppendixA}.
\end{lemma}

\begin{remark}\label{remark1} We emphasize that $(H^t)^{-1}$ is for notational purposes and no matrix inversion is needed in all cases (the exact computation scheme will be specified in the subsequent subsections). Note that to satisfy the requirement of Lemma \ref{lemma1}, zero initialization for all variables suffices. Lemma \ref{lemma1} establishes that updates (\ref{lemma1_1})--(\ref{lemma1_4}) are equivalent to (\ref{primal})--(\ref{zupdate}) and (\ref{admm4})--(\ref{admm5}) under appropriate initialization. This has a twofold implication: (i) we have achieved transforming a 3-block ADMM to a 2-block ADMM, which allows for a broader range of algorithm parameters $\mu_z,\mu_\theta$ that guarantee convergence; (ii) it is not required to explicitly store and update $z^{t}$ since it evolves on a linear manifold parameterized by $x^t$, i.e., $z=\tfrac{1}{2}E_ux$. Besides, only half of $y^t$ needs to be stored since $y^t = [\alpha^t;-\alpha^t]$. This further reduces associated storage and communication costs. We note that a 2-block ADMM can be achieved directly without introducing $z$ variables. However, such a direct formulation induces additional communication rounds when curvature information is computed. We further discuss this in \textit{Remark} 2.
\end{remark} 
Using the equivalent while more efficient updates (\ref{lemma1_1})--(\ref{lemma1_4}), we proceed to develop a family of algorithms by explicating different choices of $J^t$ used in the construction of the approximated Hessian of the augmented Lagrangian as:
\begin{gather}
    H^t = J^t+\mz D+\mt S S^\top+\epsilon I_{md}, \label{H_def}
\end{gather}
where we have introduced $\epsilon>0$ to provide additional robustness for our approximation. Notice that $\mz D +\mt S S^\top+\epsilon I_{md}$ is a diagonal matrix, whence $H^t$ is block-diagonal when $J^t$ is. When $H^t$ is \emph{block-diagonal}, each component of the update direction, $(H^t)^{-1}\nabla_x \mathcal{L}^t$ in (\ref{primal}) and equivalently in (\ref{lemma1_1}), can be computed individually by agents. More precisely, agent $i$ computes the update $u_i^t$ by solving the following linear system:
\begin{gather}
    H^t_{ii}u_i^t=\nabla_x \mathcal{L}^t_i. \label{linear_system}
\end{gather}
Therefore, once the right-hand side of (\ref{linear_system}) is obtained by $i$-th agent, no additional communication is needed to solve for $u_i^t$. This is made possible by using the intermediate consensus variables $\{z_{ij}\}$ which decouple $x_i$ from $x_j$.

\begin{remark}\label{remark2}
If consensus constraints are enforced directly as $x_i=x_j$, e.g., $E_s x=0$, then the Hessian of the augmented Lagrangian will not be block-diagonal, but rather have a structure compatible with the graph: 
\begin{gather*}
    H^t=\nabla^2 F(x^t)+\mz L_s+\mt SS^\top+\epsilon I_{md}.
\end{gather*}
Due to the presence of the signed graph Laplacian matrix $L_s$, the $ij$-th block will be nonzero if $(i,j)\in \mathcal{E}$. In such scenarios, computing $u_i^t$ either requires the presence of a fusion center that gathers all $\nabla_x\mathcal{L}^t_i$ for centralized processing, or a distributed implementation can be pursued by computing inexact (quasi) Newton-updates by truncating the Taylor expansion of the Hessian inverse with $K$ terms \cite{2nd2017,dnadmm, 2nd2020,2ndnonsmooth2020,Mokhtari2016,Eisen2017,Eisen2019,bajovic2017}. However, the truncation approach incurs $K$ additional communication rounds among agents and their neighbors, per iteration. This not only induces large communication overhead, but also demands stringent synchronization among agents \cite{Freris2011}. In contrast, all our proposed methods feature minimal communication complexity (see step \ref{comm_algo} of Alg. 1), and are amenable to an asynchronous implementation. 
Different choices of $J^t$ in (\ref{H_def}) affect the local computational cost and convergence rate as we elaborate next. 
\end{remark}

\subsection{Gradient updates}
By choosing $J^t_{\mathrm{Gradient}}\equiv0$, it follows that $H^t$ is diagonal. Therefore, (\ref{lemma1_1}) is equivalent to performing diagonally preconditioned gradient descent on the augmented Lagrangian, where step sizes are controlled by setting $\epsilon$. We note that the proposed algorithm recovers Decentralized Linearized ADMM (DLM) \cite{Ling2015} with $g(\cdot)=0$ as a special case of (\ref{prob3}). Specifically, agent $i$ computes  $u_i^t$ from (\ref{linear_system}) as:
\begin{gather}
   u_i^t =  (\mz\abs{\mathcal{N}_i}+\delta_{il}\mt  +\epsilon)^{-1} \nabla_x \mathcal{L}^t_i, \label{gradient_step_size}
\end{gather}
where $\delta_{il}=1$ if $i=l$ and $0$ otherwise. The above shows that the step size of the gradient descent at agent $i$ is related to the number of its neighbors and can be adjusted by tuning $\epsilon$. Computing updates for agents using (\ref{gradient_step_size}) involves $\mathcal{O}(d)$ computational costs, and we proceed to specify how curvature information is incorporated with nonzero $J^t$ in the following. 
\subsection{Newton updates}
By setting $J^t_{\mathrm{Newton}}=\nabla^2 F(x^t)$ in (\ref{H_def}), we obtain the Hessian of the augmented Lagrangian plus $\epsilon I_{md}$:
\begin{gather}
    H^t_{\text{Newton}} = \nabla^2 F(x^t)+\mu_z D+\mt S S^\top+\epsilon I_{md}.\label{hessian}
\end{gather}
We note that since $F(x^t)=\sum_{i=1}^m f_i(x_i^t)$, $\nabla^2 F(x^t)$ is a block diagonal matrix with the $i$-th block being $\nabla^2 f_i(x_i^t)$. As discussed previously, this induces a \emph{block diagonal} $H^t$ and the update direction $u_i^t$ can be obtained by solving (\ref{linear_system}) by each agent at the cost of $\mathcal{O}(d^3)$ for general objective functions, \emph{without additional communication among agents.}

\subsection{Quasi-Newton updates}
\hide{We further exploit the \emph{block-diagonal} structure of the Hessian by proposing the following quasi-Newton updates that remove the need for solving linear systems (\ref{linear_system}). As mentioned in Section \ref{intro_quasi}, a direct application of BFGS does not allow distributed implementation as global operations are required to construct the curvature model. Efforts were devoted in the literature to devise distributed BFGS variants \cite{Eisen2017,Eisen2019,bajovic2017}. However, these methods involve additional communication rounds among agents.} 
In this section, we introduce a distributedly implementable BFGS scheme that harnesses curvature information without inner communication loops. Some insights can be gained by investigating the target Hessian of the augmented Lagrangian in (\ref{hessian}). We note that the only time-varying part of $H^t_\mathrm{Newton}$ is $\nabla^2 F(x^t)$, while the remaining part is constant (the graph structure is assumed to be time-invariant in this paper). In \cite{Eisen2019}, authors propose to estimate $\nabla^2 F(x^t)$ using the BFGS formula with each node's local information and then compute the $K$-th order Taylor expansion of $(H^t)^{-1}$. For a distributed implementation, $K$ additional communication rounds are needed due to direct coupling between $x_i$ and $x_j$ in \cite{Eisen2019}. We note that such schemes not only incur higher communication costs per round, but also induce $\mathcal{O}(d^3)$ computational costs since linear systems have to be solved by agents. 

In contrast, we exploit the block-diagonal structure of the Hessian (\ref{hessian}), and propose the following scheme for approximating $(H^t)^{-1}$ using no additional communication (i.e., by means of local computation with information already available at the agents). In specific, each agent $i$ constructs the Hessian inverse model directly using the pairs $\{q^t_i,s^t_i\}_{i=1}^m$ defined as: 
\begin{equation}
\begin{aligned}
q^t_i := \nabla f_i(x^{t+1}_i)&-\nabla f_i(x^{t}_i)\\+
\left(\mz\abs{\mathcal{N}_i}+\delta_{il}\mt+\epsilon\right)s_i^t,&\,\,\text{and}\,\,
    s_i^t:=x_i^{t+1}-x_i^t.\label{diff}
\end{aligned}
\end{equation}
In other words, instead of approximating the Hessian inverse of the local objective $\left(\nabla^2 f_i(x_i^t)\right)^{-1}$, we are directly constructing a model for $\left(\nabla^2_{x}\mathcal{L}_i^t\right)^{-1}$. The $i$-th block of the approximated Hessian inverse $(H^{t}_{ii})^{-1}$ can be recursively updated using (\ref{bfgs_formula}) for the $\{q,s\}$ pairs defined in (\ref{diff}). We emphasize that it is not needed to explicitly form $H^t_{ii}$ and solve for the update direction as in (\ref{linear_system}). Instead, computing $u_i^t$ is tantamount to performing matrix multiplication $(H^t_{ii})^{-1} \nabla_x \mathcal{L}^t_i$. In summary, the proposed algorithm is advantageous compared to existing methods over the following aspects: (i) no additional communication loops are needed after each gradient evaluation and (ii) the computation costs for each agent is reduced from $\mathcal{O}(d^3)$ to $\mathcal{O}(d^2)$. For the sake of comparison with the gradient and Newton updates, we define:
\begin{gather}
    J^t_{\mathrm{BFGS}}:=H_{\mathrm{BFGS}}^t-\mz D-\mt S S^\top-\epsilon I_{md},\label{J_def_bfgs}
\end{gather}
where $H^t_{\mathrm{BFGS}}$ is obtained by the BFGS formula with $\{q_i^t,s_i^t\}_{i=1}^m$ pairs defined in (\ref{diff}). We proceed to describe the distributed implementation of the proposed algorithms.

\section{Asynchronous description}\label{section4}
In synchronous algorithms, all agents communicate with their neighbors and participate in computing in a coordinated and deterministic fashion. Such settings are appropriate when abundant communication bandwidth is available and the network is homogeneous in the sense that agents are able to finish local computations in adjacent time windows. In heterogeneous networks, where agents have different hardware conditions and different volumes of data, the progress of synchronous algorithms is limited by the slowest agent in the network at each iteration (also known as the straggler problem). Moreover, the requirement of a central coordinator becomes less practical when the size of the network grows and the availability of agents becomes unpredictable. 

Asynchronous algorithms \cite{Bert1989} remove the need for a central clock by letting a subset of agents update in a randomized fashion at each iteration. Asynchronous methods can be further classified into totally asynchronous algorithms and partially asynchronous. In the former setting, agents are able to tolerate arbitrarily large delays between updates while in the latter, a maximum delay constraint is imposed to guarantee convergence. In this section, we extend DRUID to the \textit{totally asynchronous} setting that further broadens its applicability. 

Recall the synchronous updates defined in (\ref{updates}). With any choice of computing scheme (gradient descent, Newton, or BFGS), we compactly express the synchronous algorithm by defining the operator $T:\mathbb{R}^{(2m+2)d}\to \mathbb{R}^{(2m+2)d}$ as follows:
\begin{gather}
    v^{t+1}= Tv^t, \label{operator}
\end{gather}
where $v\in \mathbb{R}^{(2m+2)d}$ is a concatenation of $[x;\phi;\theta;\lambda]$, and the operator $T$ maps $[x^t;\phi^t;\theta^t;\lambda^t]$ to $[x^{t+1};\phi^{t+1};\theta^{t+1};\lambda^{t+1}]$ according to (\ref{updates}). We proceed to define the following activation matrix:
\begin{gather}
    \Omega^{t}: =\begin{bmatrix}
    X^{t}   & 0 & 0&0\\
    0       & X^t & 0 & 0 \\
    0       & 0 & X^t_{ll} & 0 \\
    0       & 0 & 0 & X^t_{ll}
    \end{bmatrix},\label{activation_matrix}
\end{gather} 
where $X^t\in \mathbb{R}^{md\times md}$ is a diagonal random matrix with sub-blocks $X^t_{ii}\in\mathbb{R}^{d\times d}$, $i\in[m]$, being random sub-matrices corresponding to the $i$-th agent and taking values as the identity matrix $I_d$ or a zero matrix. Using the definition (\ref{operator}) and (\ref{activation_matrix}), the proposed asynchronous algorithms are expressed as:
\begin{gather}
    v^{t+1}= v^t+\Omega^{t+1}(Tv^t-v^t).\label{asy_algo}
\end{gather}

The above construction corresponds to activating agents, i.e., the $i$-th agent only updates the corresponding pair $(x^t_i,\phi^t_i)$ (additionally $(\theta^t,\lambda^t)$ if $i=l$) if and only if $X^{t+1}_{ii}=I_d$. We proceed to describe the implementation details of DRUID.

\begin{algorithm}[t]
	\caption{DRUID} 
	\textbf{Initialization}: zero initialization for all variables.
	\begin{algorithmic}[1]
		\For {$t=0,1,\dots$}
		\For {all active agents $i$} \label{loop_start_algo}
		\Statex\textit{Compute the local curvature $H_{ii}^t$}:\tikzmark{top}
		\State \label{cur_est_start}
          $J_{ii}^t=
		\begin{cases}
		0 & \text{Gradient updates}\\
		\nabla^2 f_i(x_i^t) & \text{Newton updates} 
		\end{cases}$ 
		\State  \label{cur_est_end}$(H^t)_{ii}\leftarrow J^t_{ii}+\left(\mu_z\abs{\mathcal{N}_i}+\delta_{il}\mt+\epsilon\right)I_d$\tikzmark{right} \tikzmark{bottom}
		\Statex \textit{Primal update}: 
		\State \label{h_algo}Compute $h_i^t$ as in (\ref{grad_def})
		\State \label{u_i^t}$\begin{cases}
		H^t_{ii} u_i^t = h_i^t & \text{Gradient/Newton updates} 	\\
		u_i^t=(H^t_{ii})^{-1} h_i^t &\text{BFGS updates}
		\end{cases}$
		\State\label{primal_algo}$x_i^{t+1}=x_i^t-u_i^t$
		\Statex \textit{Communication}: 
		\State \label{comm_algo} Broadcast $x_i^{t+1}$ to neighbors
		\Statex \textit{Dual update}:
		\State \label{dual_algo} $\phi_i^{t+1}= \phi_i^t+\tfrac{\mu_z}{2}\sum_{j\in \mathcal{N}_i}(x_i^{t+1}-x_j^{t+1})$
		\Statex \textit{Updates pertaining to the regularization function}:
		\If{$i=l$}
		\State \label{proximal_start} $\theta^{t+1}=\textbf{prox}_{g/\mt}(x_l^{t+1}+\tfrac{1}{\mt}\lambda^t)$
		\State \label{proximal_end}$\lambda^{t+1}=\lambda^t+\mt(x^{t+1}_l-\theta^{t+1})$
		\EndIf 
		\Statex \textit{Curvature estimation update (BFGS only)}:
		\State \label{bfgs_update}Update $(H^{t+1}_{ii})^{-1}$ using $\{q_i^t,s_i^t\}$ in (\ref{diff}) and the formula in (\ref{bfgs_formula}) 
		\EndFor
		\EndFor 
		\AddNote{top}{bottom}{right}{only for gradient or\\ Newton \\updates}
	\end{algorithmic} \label{algo_sync}
\end{algorithm}

\subsection{Distributed  and Asynchronous Implementation} 
 
The proposed algorithms admit the exact same implementation with variable computing choices corresponding to the selection of $J^t$ in (\ref{H_def}), so as to incorporate curvature information or not. The unified description is detailed in Algorithm 1. \hide{We note that the dual variable $\alpha^t\in \mathbb{R}^{nd}$ is only invoked in the form of $E_s^\top \alpha^t$ in the primal update (\ref{lemma1_1}). By pre-multiplying both sides of (\ref{lemma1_3}) with $E_s^\top$ and defining $\phi^t := E_s^\top \alpha^t\in \mathbb{R}^{md}$, we obtain the following update rule (recall $L_s=E_s^\top E_s$):
\begin{gather}
    \phi^{t+1}= \phi^t+\tfrac{\mz}{2}L_sx^{t+1}.
\end{gather}
This further eliminates the need for communicating dual variables by letting each agent hold the pair}We let the $i$-th agent hold $(x_i^t,\phi_i^t)$ while the $l$-th agent additionally holds the pair $(\theta^t,\lambda^t)$ pertaining to the nonsmooth regularization function $g(\cdot)$. The gradient of the augmented Lagrangian pertaining to agent $i$, $\nabla_{x}\mathcal{L}_i^t$, is expressed as: 
\begin{equation}
\begin{aligned}
    h_i^t = \nabla  f_i(x_i^t)&+\phi_i^t+\tfrac{\mu_z}{2}\sum_{j\in\mathcal{N}_i}(x_i^t-x_j^t)\\
    +\delta_{il}\mt&(x_i^t-\theta^t+\tfrac{1}{\mt}\lambda^t).
\end{aligned}\label{grad_def}
\end{equation}
Before we present the asynchronous implementation (Alg. 1), we describe the synchronous case as a special case to shed some light on the design principles. At the beginning of each round, all agents become active and estimate their local curvatures as in steps \ref{cur_est_start}-\ref{cur_est_end} (without communication, irrespective of the computing schemes). For the BFGS computing scheme, no computation is required at these steps. Agents then carry primal updates by first computing $h_i^t$ expressed in (\ref{grad_def}). We emphasize that $h_i^t$ can be computed without communication, since each agent $i$ already has access to the variables of its neighbors, $\{x_j^t |j\in\mathcal{N}_i\}$, from the previous round with zero initialization. \hide{It thus computes $h_i^t$ locally in step \ref{h_algo}.} If gradient or Newton updates is opted as the computing scheme, agents compute $u_i^t$ by solving the linear system (for gradient descent $u_i^t$ can be trivially solved since $H^t_{ii}$ is a constant scalar times the identity matrix). For the BFGS scheme, $u_i^t$ is computed by performing matrix-vector multiplication $(H^t_{ii})^{-1}h_i^t$. Once $u_i^t$ is obtained, agents update their $x_i^{t+1}$ in step \ref{primal_algo}. We note that the \textit{only communication round} occurs at step \ref{comm_algo} where agents broadcast $x_i^{t+1}$ to their neighbors (thus incurring the same cost for all computing schemes, i.e., $|\mathcal{N}_i|d$ for agent $i$). Dual updates are executed in step \ref{dual_algo}. We require agents to store $\{x_j^{t+1}$, $j\in\mathcal{N}_i\}$, to execute step \ref{h_algo} in the next iteration. In addition to the primal-dual variables $(x_l,\phi_l)$, the $l$-th agent further holds $(\theta,\lambda)$ associated with the regularization term $g(\cdot)$, which are updated in steps \ref{proximal_start}-\ref{proximal_end}. Finally, if BFGS is opted as the updating scheme, agents update \hide{their}local curvature estimation $(H_{ii}^{t+1})^{-1}$ in step \ref{bfgs_update}. 

In the case of asynchronous implementation, we equip each agent with a buffer so that even if agents are not active, they can still receive information from their neighbors. \hide{We let each agent hold a local clock, which ticks according to a Poisson process, independent of each other.\hide{We assume that clocks tick at a slower time scale than agents' computing processes, and only one clock ticks at each iteration.} These assumptions are standard in asynchronous algorithms, see also in \cite{asyn2011,poisson2015,arock2016,2nd2020}. In this setting, the $i$-th agent becomes active only if the corresponding clock ticks, captured by the activation matrix $X^{t+1}_{ii} = I_d$.} Once active, the $i$-th agent executes steps \ref{cur_est_start}-\ref{cur_est_end} using only local information and then retrieves the most recent $x^t_j$ from its buffer for computing $h_i^t$ in step \ref{h_algo}. Once $u_i^t$ is computed and $x_i^{t+1}$ is updated in steps \ref{u_i^t} and \ref{primal_algo}, respectively, the active agent $i$ broadcasts $x_i^{t+1}$ to its neighbors, whose buffers store the updated $x_i^{t+1}$.\hide{Since only one agent is active during each iteration, $x_j^{t+1}=x_j^t$ for all $j\neq i$.} Finally, the active agents check their buffers for most recent $x_j^{t+1}$ and proceed to dual updates and finish their computing as in steps \ref{dual_algo}-\ref{bfgs_update}.

\hide{\begin{itemize}
    \item Synchronous algorithms: By setting $X_{ii}^{t}=Y_{jj}^{t}=I_d$ for all $i\in [m]$, $j\in[n]$, and $t\geq 0$, we obtain $\Omega^{t}=I$ for all $t$. Therefore, the update (\ref{asy_algo}) reduces to (\ref{operator}), which corresponds to the case when all agents update in a synchronous fashion. 
    \item Single agent activation: At each iteration, the random matrix $X^{t}$ and $Y^t$ are of the form $\textbf{diag}\{\Phi_m^{t}\otimes I_d\}$ and $\textbf{diag}\{\Phi_n^{t}\otimes I_d\}$ respectively, where $\Phi_m^{t}\in\mathbb{R}^m$ takes values in the sample space: \{$\underbrace{(1,0,\dots,0)}_{m\,\,\text{entries}}$,\,(0,1,0,\dots,0),\,(0,\dots,0,1)\}, similarly for $\Phi^t_n\in\mathbb{R}^n$ with $n$ entries. 
    \item Multiagent activation: For $t\geq 0$, $i\in[m],j\in[n]$, $X_{ii}^{t}=Y_{jj}^{t}=I_d$ with probability $p_i^X,p_j^Y>0$, independent of other $i\in[m]$ and $j\in[n]$. We assume this setting of activation for our analysis in Section \ref{section5}. 
\end{itemize}}

\section{Analysis}\label{section5}
In this section, we present a unified framework for analyzing the proposed algorithms with gradient, Newton, and BFGS updates. Throughout this section, we assume that the initialization requirement in Lemma \ref{lemma1} is satisfied. We recall the concatenated vector $v=[x;\phi;\theta;\lambda]\in\mathbb{R}^{(2m+2)d}$ introduced in (\ref{operator}), and we similarly define $v_{\alpha}=[x;z;\alpha;\theta;\lambda]\in\mathbb{R}^{(m+2n+2)d}$. We use $v$ for implementation as in Algorithm 1 but analyze convergence using $v_\alpha$ for technical convenience. We note that 
their equivalence is established by Lemma \ref{lemma1} using $\phi=E_s^\top\alpha, z=\tfrac{1}{2}E_ux$. We first establish the sublinear convergence rate of the synchronous DRUID under the assumption that the local objective functions are convex. By further assuming strong convexity, we establish the global linear convergence rate for both the synchronous and the asynchronous settings. 

\subsection{Preliminaries}
\hide{We begin by stating the assumptions for our analysis.} 
\begin{assumption}\label{assumption_existence}(Existence of solutions) The solution set $\mathcal{X}^\star$ of problem (\ref{prob1}) is nonempty, i.e., $\mathcal{X}^\star\neq \emptyset$. 
\end{assumption}
\begin{assumption}\label{assumption1} The local costs functions $f_i(\cdot)$ and the regularizer function satisfy the following conditions:\\
(i) Each $f_i(\cdot):\mathbb{R}^d\to \mathbb{R}$ is twice continuously differentiable, $m_f$--strongly convex and $M_f$--smooth, i.e., $\forall i\in[m]\,,x_i\in \mathbb{R}^d$: 
\begin{gather}
    m_f I_d \preceq \nabla^2 f_i(x_i) \preceq M_f I_d, \label{f assumption}
\end{gather}
where $0\leq m_f\leq M_f<\infty$.\\
(ii) The regularizer function $g(\cdot):\mathbb{R}^d\to \mathbb{R}$ is proper, closed, and convex, i.e., $\forall\,x,y\in\mathbb{R}^d$,
\begin{gather}
    (x-y)^\top (\partial g(x)-\partial g(y))\geq 0, \label{g assumption}
\end{gather}
where the inequality is meant for arbitrary elements in the subdifferential sets $\partial g(x)$ and $\partial g(y)$, respectively. 
\end{assumption}

\begin{assumption}\label{lip_hessian}
The Hessians of the local objective functions are Lipschitz continuous with constant $L_f$, i.e., $\forall\,i\in[m],x,y\in \mathbb{R}^d$,
\begin{align*}
    \norm{\nabla^2 f_i(x)-\nabla^2 f_i(y)}\leq L_f\norm{x-y}.
\end{align*}
\end{assumption}
\noindent Note that we allow the case $m_f=0$ (convex but not strongly convex), and we will analyze separately for the cases $m_f=0$ and $m_f>0$ to establish sublinear and linear rates, respectively. \hide{Moreover, since $F(x)=\sum_{i=1}^m f_i(x_i)$, the lower and upper bound in (\ref{f assumption}) hold for $\nabla^2 F(x)$ (with $I_d$ replaced with $I_{md}$), as well.}Assumptions \ref{assumption_existence}-\ref{assumption1} are standard for analyzing distributed algorithms while Assumption \ref{lip_hessian} is standard for analyzing second-order methods \cite{Nesterov2018}.  

\begin{assumption}\label{assumption_bounded_bfgs}
The Hessian estimate obtained by the BFGS is uniformly upper bounded, i.e., for any $t\geq0$, there exists a constant $\psi>0$ such that:
\begin{align}
    H^t_{\mathrm{BFGS}} \preceq \psi I_{md}. \label{bfgs_bound}
\end{align}
\end{assumption}

\begin{remark}
Assumption \ref{assumption_bounded_bfgs} applies only for BFGS updates and is, in general, not standard. However, many techniques can be used to satisfy (\ref{bfgs_bound}). For example, adding small regularization when computing the Hessian inverse approximations, i.e., $(H^t_{\mathrm{BFGS}})^{-1}=(\widehat{H}^t_{\mathrm{BFGS}})^{-1}+\tfrac{1}{\psi}I_{md}$, where $(\widehat{H}^t_{\mathrm{BFGS}})^{-1}$ is obtained through (\ref{bfgs_formula}). Other means\hide{for ascertaining the upper bound in (\ref{bfgs_bound})} include using regularized BFGS updates and invoking L-BFGS\cite{icassp2022} estimation by using a finite prescribed number of $\{q_i^t,s_i^t\}$ copies. In brief, we make this assumption for convenience and without serious loss in generality; see also \cite{Eisen2019} and \cite{Mokhtari2014}. 
\end{remark} 

When local functions are assumed to be only convex ($m_f=0$), there might be multiple optimal primal solutions\hide{$x^\star$}, each with multiple optimal dual solution\hide{$\alpha^\star$}.\hide{For each optimal primal solution $x^\star$, there are multiple optimal dual solutions $\alpha^\star$.} However, there exists a unique dual pair that lies in the column space of some matrix, to be defined and formalized in the following. \hide{We formalize this in the following lemma.}  
\begin{lemma}\label{optimal}
The tuple $(x^\star,z^\star,\alpha^\star,\theta^\star,\lambda^\star)$ solves (\ref{prob3}), and equivalently (\ref{prob1}), if and only if the following holds:
\begin{align}
    \nabla F(x^\star)+E_s^\top \alpha^\star+S \lambda^\star &= 0, \tag*{KKTa}\label{KKTa}\\
    \partial g(\theta^\star)-\lambda^\star &\ni 0 ,\tag*{KKTb}\label{KKTb}\\
    E_s x^\star& =0,\tag*{KKTc}\label{KKTc}\\
    E_u x^\star& = 2z^\star ,\tag*{KKTd}\label{KKTd}\\ 
    S^\top x^\star&=\theta^\star.\tag*{KKTe}\label{KKTe}
\end{align}
Moreover, there exists a unique dual optimal pair $[\alpha^\star;\lambda^\star]\in \mathbb{R}^{(n+1)d}$ that lies in the column space of $C:=\begin{bmatrix}
E_s \\ S^\top
\end{bmatrix}\in \mathbb{R}^{(n+1)d\times md}$.

\textit{Proof}: See Appendix \ref{AppendixA} of the extended version of this paper.
\end{lemma}

\hide{The above conditions can be used to characterize the suboptimality of iterates $(x^t,z^t,\alpha^t,\theta^t,\lambda^t)$ generated by (\ref{updates}), where KKTa-b can be interpreted as the optimality condition of minimizing the objective function, while KKTc-e measure consensus violation among agents. }We proceed to establish a lemma that characterizes the suboptimality of the iterates when replacing (\ref{admm1}) with\hide{one-step updates as in} (\ref{lemma1_1}). \hide{The following lemma also serves to relate primal iterates $(x^{t+1},\theta^{t+1})$ to dual iterates $(\alpha^{t+1},\lambda^{t+1})$.} 
\begin{lemma}\label{primal_sub}Consider the iterates generated by (\ref{updates}). The following holds: 
\begin{gather}
    e^t+\nabla F(x^{t+1})-\nabla F(x^\star)+\epsilon (x^{t+1}-x^t)+E_s^\top (\alpha^{t+1}-\alpha^\star)\nonumber\\
    +\mz E_u^\top(z^{t+1}-z^t)
   +S\left( \lambda^{t+1}-\lambda^\star+\mt\left( \theta^{t+1}-\theta^t\right)\right) = 0\nonumber
\end{gather}
where the error term is:
\begin{gather}
        e^t = \nabla F(x^t)-\nabla F(x^{t+1})
        +J^t(x^{t+1}-x^t). \label{error}
\end{gather}
\textit{Proof}: See Appendix \ref{AppendixA}.
\end{lemma}
\hide{We proceed to establish convergence of DRUID and characterize the rate in the next subsection.}
\hide{It is not hard to verify that Lemma \ref{primal_sub} holds with $e^t=0$ when exact optimization is performed as in (\ref{admm1}). Since the difference between the proposed algorithms lies in how the primal update (\ref{lemma1_1}) is executed with different choices of curvature estimation, characterizing how the error in (\ref{error}) behaves for variable methods is the key to reveal the differences between algorithms with gradient, Newton, and BFGS updates. We delay the discussion on how $e^t$ varies with different updating schemes to the next subsection.} \hide{The following lemma shows that KKTb holds for any pair $(\theta^{t+1},\lambda^{t+1})$ generated by (\ref{updates}).} 

\hide{\begin{lemma}\label{inclusion lemma}
Consider the dual update $\theta^{t+1}=\textbf{prox}_{g/\mu_\theta}( S^\top x^{t+1}+\tfrac{1}{\mu_\theta}\lambda^t)$ (\ref{lemma1_2}) and any optimal pair $(\theta^\star,\lambda^\star)$. The following holds:
\begin{align}
\partial g(\theta^{t+1})-\lambda^{t+1}\ni 0, &\\
 (\lambda^{t+1}-\lambda^\star)^\top (\theta^{t+1}-\theta^\star) \geq 0,&\label{claim1_1}\\
 (\lambda^{t+1}-\lambda^t)^\top (\theta^{t+1}-\theta^t)\geq 0.\label{claim1_2}&
\end{align}
\textit{Proof}: See Appendix \ref{AppendixA}.
\end{lemma}
Since KKTb holds along the iterates path as shown by Lemma \ref{inclusion lemma}, it suffices to show that the optimality residual of KKTa,c,d converge to $0$ to establish convergence.}

\subsection{Sublinear Convergence}
We recall $J^t$ in (\ref{H_def}) and the concatenated vector $\va\in\mathbb{R}^{(m+2n+2)d}$. We further define $\overline{J}^t = J^t+\epsilon I$,  and the scaling matrix $\mathcal{G}^t$ as follows:
\begin{gather}
    v_\alpha = \begin{bmatrix}
    x\\ z\\ \alpha \\ \theta \\ \lambda
    \end{bmatrix}, \quad \mathcal{G}^t = 
    \begin{bmatrix}
    \,\overline{J}^t & 0    & 0             & 0   & 0\\
    0                & 2\mz & 0             & 0   & 0\\
    0                & 0    &\tfrac{2}{\mz} & 0   & 0 \\
    0                & 0    & 0             & \mt & 0\\
    0                & 0    & 0             &0    &\tfrac{1}{\mt}
    \end{bmatrix}. \label{u_g_def}
\end{gather}
\begin{theorem}\label{sublinear_theorem_combined} Recall the definition in (\ref{u_g_def}). Consider the iterates generated by (\ref{updates}). We denote the smallest and the biggest eigenvalue of $L_u$ and $L_s$ as $\sigma_{\mathrm{min}}^{L_u}$ and $\sigmaxLs$ respectively. Under Assumptions \ref{assumption_existence}-\ref{assumption_bounded_bfgs}, ($m_f=0$), and we select $\mz$ and $\epsilon$ such that:
$
    \epsilon>\tfrac{M_f}{2}, \mz\epsilon<\psi^2.
$
Then the following holds:
\begin{align}
      &\tfrac{1}{T}\tfrac{\mz}{2}\norm{x^1}^2_{L_s}+\tfrac{\mt}{T}\norm{x^1_l-\theta^1}^2 +\tfrac{1}{T}\sum_{t=1}^T\norm{\va^{t+1}-\va^t}^2_{\mathcal{G}^t}\label{sublinear_theorem_2}\\
      &\geq  \tfrac{1}{T}\tfrac{\mz}{2}\norm{x^{T+1}}^2_{L_s}+\tfrac{\mt}{T}\norm{x^{T+1}_l-\theta^T}^2\nonumber\\
      &+\tfrac{1}{T}\sum_{t=1}^T \Big\{\tfrac{\epsilon}{\rho\overline{M}^2}\norm{\nabla F(x^t)+E_s^\top \alpha^t+S \lambda^t}^2+\mt\norm{\theta^{t+1}-\theta^t}^2\nonumber\\
      &+2\mz\norm{z^{t+1}-z^t}^2
      +\left(\tfrac{\mz}{2}-\tfrac{\epsilon\mz^2}{2\overline{M}^2}\right)\norm{x^t}^2_{L_s}\nonumber\\
      &+\left(\mt-\tfrac{2\epsilon\mt^2}{\overline{M}^2(\rho-1)}\right)\norm{ S^\top x^t-\theta^t}^2\Big\}, \nonumber
\end{align}
where\hide{we denote} $d_{\mathrm{max}}=\underset{i}{\mathrm{max}}\,\, \abs{\mathcal{N}_i}$, $\rho> \max\left\{\tfrac{2\epsilon \mt}{\overline{M}^2},\sigmaxLs\right\}+1$, and $\overline{M}$ (for each scheme) is given by:
$
    \overline{M}_{\mathrm{Gradient}} =\mz d_\mathrm{max}+\epsilon+\mt, 
    \overline{M}_{\mathrm{Newton}} =M_f+\mz d_\mathrm{max} +\epsilon +\mt,
    \overline{M}_{\mathrm{BFGS}} =\psi.
$ 

\textit{Proof}: See Appendix \ref{AppendixB}.
\end{theorem}

\hide{\begin{theorem}\label{descent_theorem} Consider the iterates generated by (\ref{updates}). Denote the smallest and the biggest eigenvalue of $L_u$ and $L_s$ as $\sigma_{\mathrm{min}}^{L_u}$ and $\sigmaxLs$ respectively. Under Assumptions \ref{assumption_existence}-\ref{assumption_bounded_bfgs}, (convex $f_i(\cdot)$ with $m_f=0$), and we select $\mt,\mz$, and $\epsilon$ such that
\begin{align}
    \tfrac{\mz\sigma_{\mathrm{min}}^{L_u}}{2}+\epsilon-\tfrac{M_f}{2}&>0,\\
    2(m+\tfrac{\mt}{\mz})^2&>\sigma^{L_u}_{\mathrm{min}}.
\end{align}
Then the following holds:
\begin{align}
  \norm{v^{t+1}-v^\star}^2_{\mathcal{G}^t}\leq \norm{v^t-v^\star}^2_{\mathcal{G}^t}-\delta \norm{v^{t+1}-v^t}^2_{\mathcal{G}^t},\label{sublinear_theorem_1}
\end{align}
where $\delta= \left(\tfrac{\mz\sigma_{\mathrm{min}}^{L_u}}{2}+\epsilon-\tfrac{M_f}{2}\right)/\left(\tfrac{\mz\sigma_{\mathrm{min}}^{L_u}}{2}+\epsilon\right)>0$.\\
\textit{Proof}: See Appendix \ref{AppendixB}.
\end{theorem} 

Theorem \ref{descent_theorem} shows that the Lyapunov function $\norm{v^{t}-v^\star}^2_{\mathcal{G}^t}$ is monotonically decreasing and therefore converges. Standard analysis techniques can be used to establish convergence \cite{He1994} and \cite{Yin2017}. We proceed to characterize the rate of convergence. 

\begin{theorem}\label{sublinear_theorem} Recall the definition of $v$ and $\mathcal{G}^t$ defined in (\ref{u_g_def}). Under Assumptions \ref{assumption_existence}-\ref{assumption_bounded_bfgs}, the running-average suboptimality residual is upper bounded as follows: 
\begin{align}
      &\tfrac{1}{T}\sum_{t=1}^T\norm{v^{t+1}-v^t}^2_{\mathcal{G}^t}+\tfrac{1}{T}\tfrac{\mz}{2}\norm{x^1}^2_{L_s} \label{sublinear_theorem_2}\\
      &\geq  \tfrac{1}{T}\sum_{t=1}^T \big(\tfrac{\epsilon}{\rho\overline{M}_f^2}\norm{\nabla F(x^t)+E_s^\top \alpha^t+S \lambda^t}^2\nonumber\\
      &+\left(\tfrac{\mz}{2}-\tfrac{\epsilon\mz^2}{2\overline{M}^2}\right)\norm{x^t}^2_{L_s}+\left(\mt-\tfrac{2\epsilon\mt^2}{(\rho-1)\overline{M}^2}\right)\norm{S^\top x^t-\theta^t}^2\nonumber\\
    &+\mt\norm{\theta^{t+1}-\theta^t}^2\big)+\tfrac{1}{T}\tfrac{\mz}{2}\norm{x^{T+1}}^2_{L_s}+\tfrac{\mt}{T}\norm{\theta^{T+1}-\theta^T}^2, \nonumber
\end{align}
where $\rho\geq \max\{\sigmaxLs,\tfrac{\mt\mz\sigma^{L_u}_{\mathrm{min}}+2\mt\epsilon}{\overline{M}^2}\}+1$,
$
    \overline{M}_{\mathrm{Gradient}} =\mz m+\epsilon+\mt, 
    \overline{M}_{\mathrm{Newton}} =M_f+\mz m +\epsilon +\mt,
    \overline{M}_{\mathrm{BFGS}} =\psi.
$
\textit{Proof}: See Appendix \ref{AppendixB}.
\end{theorem}}

\begin{remark}
It is not hard to verify $z^t=\frac{1}{2}E_u x^t$ (\textit{Remark} \ref{remark1}) and $\lambda^t\in\partial g(\theta^t)$ holds along the convergence path and establishing convergence amounts to satisfying KKTa,c,e. \hide{The above theorem characterizes the suboptimality residuals, the consensus error, and the discrepancy between $(x^t_l,\theta^t)$, corresponding to KKTa,c,d, respectively.} We proceed to explicate the convergence rate of these terms in the following.
\end{remark}

\begin{corollary}\label{cor1}
The running-average suboptimality residual and consensus errors converge as follows:
\begin{subequations}\label{corollary1}
\begin{align}
    \tfrac{1}{T}\sum_{t=1}^T \norm{\nabla F(x^t)+E_s^\top \alpha^t+S \lambda^t}^2&= \mathcal{O}(\tfrac{1}{T}),\\
    \tfrac{1}{T}\sum_{t=1}^T\norm{x^t}_{L_s}^2&=\mathcal{O}(\tfrac{1}{T}), \\
    \tfrac{1}{T}\sum_{t=1}^T \norm{S^\top x^t-\theta^t}^2&=\mathcal{O}(\tfrac{1}{T}). 
\end{align}
\end{subequations}
\textit{Proof}: See Appendix \ref{AppendixB}. 
\end{corollary}

\subsection{Linear Convergence}

By further assuming strongly convex $f_i(\cdot)$ ($m_f>0$), we establish the linear convergence rate of DRUID. We show that the iterates converge to the unique $[x^\star;z^\star;\alpha^\star;\theta^\star;\lambda^\star]$, where the dual pair $[\alpha^\star;\lambda^\star]$ lies in the column space of $C$ as shown in Lemma \ref{optimal}. We first bound the error in (\ref{error}). 

\begin{lemma}\label{lemma_error_bound} Recall the error term defined in (\ref{error}). The following holds: 
$
    \norm{e^t}\leq \tau^t \norm{x^{t+1}-x^t},
$
where 
\begin{subequations}
\begin{align}
\tau^t_{\mathrm{Gradient}} &= M_f, \label{error_bound_gradient}\\
\tau^t_{\mathrm{Newton}} &= \min\big\{2M_f,\tfrac{L_f}{2}\norm{x^{t+1}-x^t}\big\},\label{error_bound_newton}\\
\tau^t_{\mathrm{BFGS}} & = \norm{H^{t}-H^{t+1}}\leq 2\psi.\label{error_bound_bfgs}
\end{align}
\end{subequations}
\textit{Proof}: See Appendix \ref{AppendixC}.
\end{lemma}
\noindent The above lemma complements the result presented in \cite{Ling2015} and \cite{DQM2016}. By upper bounding the error induced when we replace the exact suboptimization step (\ref{admm1}) with a one-step update (\ref{lemma1_1}), we reveal the differences when using different computing schemes. Since the algorithm converges, as established in the previous subsection, the $\tfrac{L_f}{2}\norm{x^{t+1}-x^t}$ term will eventually become smaller than the $2M_f$ term in (\ref{error_bound_newton}). In other words, the error term eventually diminishes quadratically with respect to $\norm{x^{t+1}-x^t}$ in Newton updates. On the other hand, since $\norm{H^{t}-H^{t+1}}\to 0$, we conclude that the error term in BFGS diminishes superlinearly with respect to $\norm{x^{t+1}-x^t}$. We have summarized these along with other properties of different updating schemes in Table \ref{table}. Note that the $l$-th agent, that performs updates pertaining to $g(\cdot)$, additionally holds the $(\theta, \lambda)$ pair; thus, storage increases by $2d$ and additional computation is incurred for evaluating the proximal operator (typically $\mathcal{O}(d)$). The fact that all computing schemes share the same communication cost (equal to the vector dimension $d$) is because agents only communicate once per iteration with their neighbors (step \ref{comm_algo} of Alg. \ref{algo_sync}).   \hide{In summary, as $t\to \infty$, 
\begin{subequations}\label{e_comparison}
\begin{align} 
    e^t_{\mathrm{Gradient}} &= \mathcal{O}(\norm{x^{t+1}-x^t}),\\
    e^t_{\mathrm{Newton}} &= \mathcal{O}(\norm{x^{t+1}-x^t}^2),\\
    e^t_{\mathrm{BFGS}} &=o(\norm{x^{t+1}-x^t}).
\end{align}
\end{subequations}}

\begin{table}[t]
    \caption{ Comparison between updating schemes in terms of communication and computation costs per iteration, and storage costs per agent as a function of vector dimension $d$ and neighborhood size $\abs{\mathcal{N}_i}$. The last column characterizes the decay rate of $e^t$ in Lemma \ref{lemma_error_bound} in terms of the difference $x^t_e:=x^{t+1}-x^t$.} 
    \centering
    \begin{tabular}{lcccc}
    \hline\hline
    Methods & Comm. costs & Comp. costs & Storage costs & Decay rate \\
    \hline
    Gradient  & $\abs{\mathcal{N}_i}d$ & $\mathcal{O}(d)$ & $\mathcal{O}(d)$ & $\mathcal{O}(\norm{x_e^t})$\\
    BFGS      & $\abs{\mathcal{N}_i}d$ & $\mathcal{O}(d^2)$ & $\mathcal{O}(d^2)$ & $o(\norm{x^t_e})$   \\
    Newton    & $\abs{\mathcal{N}_i}d$ & $\mathcal{O}(d^3)$ & $\mathcal{O}(d^2)$ & $\mathcal{O}(\norm{x^t_e}^2)$
    \\
    \hline
    \end{tabular}
    \label{table}
\end{table}

Before establishing the linear convergence rate of DRUID, we recall $\va$ and introduce the following diagonal scaling matrix $\mathcal{H}=\textbf{diag}[\epsilon,2\mz,\tfrac{2}{\mz},\mt,\tfrac{1}{\mt}]$ similar to (\ref{u_g_def}).
\hide{\begin{gather}\label{linear_scaling_def}
    \va=\begin{bmatrix}
    x\\ z\\ \alpha\\ \theta\\ \lambda
    \end{bmatrix},\quad\mathcal{H}= 
    \begin{bmatrix}
    \epsilon & 0    & 0             &  0 & 0\\
    0        & 2\mz & 0             &  0 & 0\\
    0        & 0    &\tfrac{2}{\mz} &  0 & 0 \\
    0        & 0    & 0             &\mt & 0 \\
    0        &  0   & 0             & 0  & \tfrac{1}{\mt}  
    \end{bmatrix}.
\end{gather}}
\begin{theorem}\label{theorem_linear_syn} Under Assumptions \ref{assumption_existence}--\ref{assumption_bounded_bfgs} with $m_f>0$, we denote the maximum and minimum eigenvalue of $L_u$ as $\sigma^{L_u}_{\mathrm{max}}$ and $\sigma^{L_u}_{\mathrm{min}}$ respectively. Let $\sigma^+_{\mathrm{min}}$ be the smallest positive eigenvalue of $C C^\top$, where $C:=\begin{bmatrix}
E_s\\ S^\top
\end{bmatrix}$, and $c_{\mathrm{max}}=2\cdot\max\{M_f,\psi\}$. By selecting $\mz=2\mt$, $\epsilon>\tfrac{c^2_{\mathrm{max}}(m_f+M_f)}{2m_fM_f}$, and arbitrary constant $\zeta\in\left(\tfrac{m_f+M_f}{2m_fM_f},\tfrac{\epsilon}{(\tau^t)^2}\right)$, the iterates generated by (\ref{updates}) satisfy:
\begin{align*}
    \norm{\va^{t+1}-\vas}^2_{\mathcal{H}} \leq \tfrac{1}{1+\eta} \norm{\va^t-\vas}^2_{\mathcal{H}},
\end{align*}
where $\eta$ satisfies:
\begin{gather}
    \eta=\min \bigg\{\left(\tfrac{2m_fM_f}{m_f+M_f}-\tfrac{1}{\zeta}\right)\tfrac{1}{\epsilon+\mt(\sigma^{L_u}_{\mathrm{max}}+2)},\tfrac{1}{2},\tfrac{2}{5}\tfrac{\mt\sigma^+_{\mathrm{min}}}{m_f+M_f},\nonumber\\
    \tfrac{\mt\sigma^+_{\mathrm{min}}(\epsilon-\zeta(\tau^t)^2)}{5((\tau^t)^2+\epsilon^2)},\tfrac{\sigma^+_{\mathrm{min}}}{5\max\{1,\sigma^{L_u}_{\max}\}}\bigg\}\label{eta}
\end{gather}
\textit{Proof}: See Appendix \ref{AppendixC}.\end{theorem}
\begin{remark}To shed some light on the convergence rate, we first consider the case when the sub-optimization problem (\ref{admm1}) is solved exactly. When an exact solution is obtained\hide{from the step (\ref{admm1})}, Lemma \ref{primal_sub} holds with $e^t=0$, and therefore $\tau^t=0$ in Lemma \ref{lemma_error_bound}. Having $\tau^t=0$ allows us to choose $\epsilon=0$ and $\zeta\gg 1$, which gives the following rate:
\begin{gather}
    \eta_{\mathrm{exact}}=\min\bigg\{\tfrac{2m_fM_f}{m_f+M_f}\tfrac{1}{\hide{\epsilon+}\mt(\sigma^{L_u}_{\mathrm{max}}+2)},\tfrac{1}{2},\tfrac{2}{5}\tfrac{\mt\sigma^+_{\mathrm{min}}}{(m_f+M_f)},\nonumber\\
    \hide{\tfrac{\mt\sigma^+_{\mathrm{min}}}{5\epsilon},}\tfrac{\sigma^+_{\mathrm{min}}}{5\max\{1,\sigma^{L_u}_{\mathrm{max}}\}}\bigg\}. \label{exact_rate}
\end{gather}
Denoting $\kappa=M_f/m_f$ and choosing $\mt = m_f\sqrt{\kappa}$, we obtain an iteration complexity of $\mathcal{O}\left(\sqrt{\kappa}\log(1/\varepsilon)\right)$ from (\ref{exact_rate}), where $\varepsilon$ is the solution accuracy and not to be confused with the hyperparameter $\epsilon$. Moreover, since $\sigma^+_{\mathrm{min}}$ is related to the smallest positive eigenvalue of $L_s$, i.e., the algebraic connectivity of the graph, (\ref{exact_rate}) implies that a more connected graph (larger $\sigma^+_{\mathrm{min}}$) gives rises to a larger $\eta_{\mathrm{exact}}$, and faster convergence rates. On the other hand, the rate $\eta$ established in (\ref{eta}) is no larger than $\eta_{\mathrm{exact}}$ in (\ref{exact_rate}): this is due to the fact that we have replaced the exact optimization step with the one step update (\ref{lemma1_1}). Characterizing the gap between $\eta$ and $\eta_{\mathrm{exact}}$ serves to reveal the differences between using gradients, Newton, and BFGS updates. This is achieved by comparing the upper bound for the error term $\tau^t$, and how $e^t$ (Lemma \ref{lemma_error_bound}) evolves as characterized by the last column of Table \ref{table}. As established in Section V-B, $\lim_{t\to\infty}\norm{x^{t+1}-x^t}=0$, the error bound $\tau^t_{\mathrm{Newton}/\mathrm{BFGS}}\to 0$ from inspecting (\ref{error_bound_newton}) and (\ref{error_bound_bfgs}). In other words, we can recover the convergence rate in (\ref{exact_rate}) only if we use curvature-aided updates.
\end{remark}

We recall that the asynchronous implementation in (\ref{asy_algo}) is defined using $\{\phi_i\}_{i=1}^m$ and $v=[x;\phi;\theta;\lambda]$, for the most efficient and economical deployment. In the rest of this section, we first characterize the condition for $v_\alpha$ to converge under random activation, and then show that the implementation (\ref{asy_algo}) satisfies this condition. We first define the following activation matrix corresponding to $v_\alpha=[x;z;\alpha;\theta;\lambda]\in\mathbb{R}^{(m+2n+2)d}$:
\begin{gather}
    \Omega_\alpha^{t}: =\begin{bmatrix}
    X^{t}   & 0   &        0 &       0 & 0\\
    0       & Y^t &        0 &       0 & 0 \\
    0       & 0   &      Y^t &       0 & 0 \\
    0       &   0 &        0 &X^t_{ll} & 0 \\
    0       &   0 &        0 & 0       &X^t_{ll}
    \end{bmatrix}.\label{activation_matrix_2}
\end{gather} The activation matrix $\Omega^t_\alpha$ differs from $\Omega^t$ in (\ref{activation_matrix}) as we allow $(z^t,\alpha^t)$ to be updated independently from $x^t$, captured by the random matrix $Y^t\in\mathbb{R}^{nd\times nd}$. We can similarly develop an asynchronous algorithm as:
\begin{gather}
    \va^{t+1}= \va^t+\Omega^t_\alpha (T_\alpha \va^t-\va^t), \label{asy_algo_alpha}
\end{gather}where the operator $T_\alpha:\mathbb{R}^{(m+2n+2)d}\to \mathbb{R}^{(m+2n+2)d}$ is equivalent to the synchronous updates (\ref{updates}). The update (\ref{asy_algo_alpha}) captures a wider range of random activation schemes than the update in (\ref{asy_algo}), but it is more costly to implement. Therefore, we only use (\ref{asy_algo_alpha}) as a guideline for analysis. We proceed to define  $\mathbb{E}^t[\cdot]:=\mathbb{E}[\cdot\vert \mathcal{F}^t]$, where $\mathcal{F}^t$ is the filtration generated by $(X^1,\hdots, X^t)$ and $(Y^1,\hdots,Y^t)$.

\begin{theorem}\label{theorem_asy} Consider the iterates generated by the asynchronous algorithm (\ref{asy_algo_alpha}). Under the same setting as in the Theorem \ref{theorem_linear_syn} and any activation scheme such that $\mathbb{E}^t[\Omega^{t+1}_\alpha]=\Omega_\alpha\succ 0$, then the following holds:
\begin{gather*}
    \mathbb{E}^t\left[\norm{\va^{t+1}-\vas}^2_{\mathcal{H}\Omega_\alpha^{-1}}\right]\leq \left(1-\tfrac{p^{\mathrm{min}}\eta}{1+\eta}\right)\norm{\va^{t}-\vas}^2_{\mathcal{H}\Omega_\alpha^{-1}},
\end{gather*}
where for $i\in[m],j\in[n]$, we denote $\mathbb{E}^t[X^{t+1}_{ii}]=p^X_i,\mathbb{E}^t[Y^{t+1}_{jj}]=p^Y_j$, $p^{\mathrm{min}}:=\underset{i\in[m],j\in[n]}{\min}\{p_i^X,p_j^Y\}$, and $\eta$ is given by (\ref{eta}). \\
\textit{Proof}: See Appendix \ref{AppendixC}.
\end{theorem}
We note that the activation of the asynchronous scheme using $(\Omega_\alpha^{t+1},\va^{t+1})$ described by (\ref{asy_algo_alpha}) amounts to specifying the random matrix $X^{t+1}$ and $Y^{t+1}$, which can be chosen independently from each other. Theorem \ref{theorem_asy} shows that as long as $\mathbb{E}^t[\Omega^{t+1}_\alpha]\succ 0$, iterates $\va^t$ converge linearly in expectation.\hide{This translates to each $x_i,i\in[m],$ and $(z_k,\alpha_k),k\in[n],$ being updated infinitely often.} On the other hand, the\hide{more practical} implementation using $(\Omega^{t+1},v^{t+1})$ described by (\ref{asy_algo}), only needs to specify $X^{t+1}$, i.e., activating agents.\hide{Note that irrespective of the activation scheme, dual updates are incorporated in the form of $\phi^{t}=E_s\alpha^{t}$.} The difference of the two\hide{activation schemes} lies in the fact that (\ref{asy_algo_alpha}) first updates a subset of $\alpha_k$, then computes $\phi=E_s^\top \alpha$, while (\ref{asy_algo}) directly updates a subset of $\phi_i$. In the following corollary, we show that using the activation scheme described by (\ref{asy_algo}), the induced iterates $v^{t}_\alpha$ converge linearly in expectation.  
\begin{corollary}\label{corollary_asy}
Consider activation matrices $\Omega^{t+1}$ in (\ref{activation_matrix}) and $\Omega^{t+1}_\alpha$ in (\ref{activation_matrix_2}). Under the same $X^{t+1}$ and updating scheme (\ref{asy_algo}), if $\mathbb{E}^t[X^{t+1}]\succ 0$, then it holds that $\mathbb{E}^t[\Omega^{t+1}_\alpha]\succ 0$. \\
\textit{Proof}: See Appendix \ref{AppendixC}. 
\end{corollary}
Since $\mathbb{E}^t[\Omega^{t+1}_\alpha]\succ 0$, we conclude that the implementation using $(\Omega^t,v^t)$ induces an equivalent sequence of $(\Omega_\alpha^t,\va^t)$ that converges linearly in expectation using Theorem \ref{theorem_asy}.

\section{Numerical experiments}

In this section, we present a comparative experimental validation of the proposed methods with existing state-of-the-art methods, namely PGE \cite{pgextra}, P2D2\cite{p2d2}, and ESOM\cite{Mokhtari2016}. Note that ESOM and other existing (quasi) Newton methods \cite{DQM2016,Eisen2017,Eisen2019,2nd2020} do not support nonsmooth regularization functions and therefore ESOM is only compared in the Fig \ref{ridge}, where the regularization function is differentiable. We consider the following distributed optimization problem:
\begin{gather}
    \underset{x\in\mathbb{R}^{d}}{\text{minimize}} \,\,G(x)=\left\{\sum_{i=1}^m f_i(x)+g(x)\right\}, \label{total_obj_function}
\end{gather}
where $f_i(\cdot)$ and $g(\cdot)$ are to be specified according to the application. All experiments are conducted using real-life data sets from the LIBSVM\footnote{https://www.csie.ntu.edu.tw/~cjlin/libsvm/} and UCI Machine Learning Repository\footnote{https://archive.ics.uci.edu/ml/index.php}. We generate connected random graphs with $m$ agents by repetitively drawing edges between agents according to a Bernoulli($p$) distribution. We ensure connectedness by redrawing the graph if necessary. The mixing matrices of P2D2 and ESOM are generated using the Metropolis rule\cite{p2d2} while the mixing matrix of PG-EXTRA is generated by the Laplacian-based constant weight matrix\cite{pgextra}, respectively. 

\subsection{Distributed LASSO}

The distributed LASSO problem considers solving (\ref{total_obj_function}) with $g(x)=\gamma \norm{x}_1$, $\gamma\in \mathbb{R}$, and $f_i(\cdot):\mathbb{R}^d\to \mathbb{R}$ defined as:
\begin{gather}
        f_i(x)=\tfrac{1}{2}\sum_{i=1}^{m_i} \left(a_i^\top x-b_i\right)^2. \label{lasso_f}
\end{gather}
Each $\{a_i,b_i\}\in\mathbb{R}^d\times \mathbb{R}$ is a given data point and $m_i$ denotes the total number of data points held by the $i$-th agent. The purpose of the regularization function $\gamma\norm{x}_1$ is to promote a sparse solution vector. We consider the Combined Cycle Power Plant (CCPP) dataset from the UCI Machine Learning Repository, using 9,000 data points of dimension $d=4$.\hide{collected from a combined cycle power plant over 6 years. The goal is to predict the net hourly electrical energy output of the plant given 4 features of the plant. We evenly distribute the data points among $m=20$ agents in a randomly generated graph with $p=0.2$ and $\gamma=0.002$.}   
\begin{figure}[t]
  \centering
  \begin{subfigure}[b]{0.49\linewidth}
    \includegraphics[width=\textwidth]{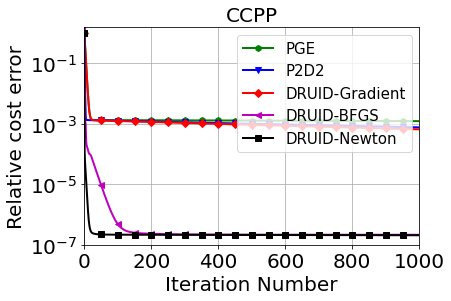}
  \end{subfigure}
  \begin{subfigure}[b]{0.49\linewidth}
    \includegraphics[width=\textwidth]{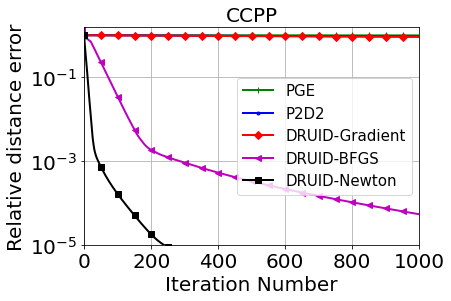}
  \end{subfigure}
  \caption{Performance comparison on the CCPP dataset. We plot the iteration number versus the relative cost error (left) $\tfrac{G(x^t)-G(x^\star)}{G(x^0)-G(x^\star)}$ and the relative distance error (right) $\tfrac{\norm{x^t-x^\star}}{\norm{x^0-x^\star}}$ on a randomly generated graph consisting of $m=20$ agents.}
  \label{lasso}
 \end{figure}
 \begin{figure}[t]
  \centering
  \begin{subfigure}[b]{0.49\linewidth}
    \includegraphics[width=\textwidth]{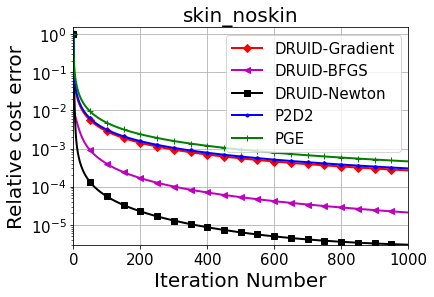}
  \end{subfigure}
  \begin{subfigure}[b]{0.49\linewidth}
    \includegraphics[width=\textwidth]{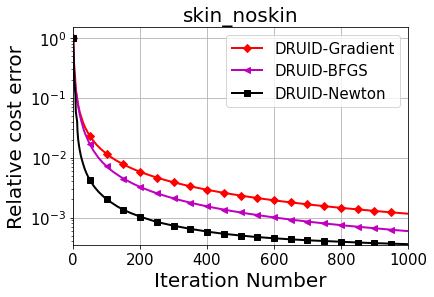}
  \end{subfigure}
  \caption{Performance comparison of DRUID algorithms and existing methods in a network with $m=20$ agents, in synchronous (left) and asynchronous (right) settings. In each iteration of the asynchronous setting, half of the agents in the network are activated in a uniformly random fashion.}
  \label{sync_async}
 \end{figure}
\indent We note that all algorithms have the same communication costs per iteration (this is due to the fact that only one round of communication of the local variable $x_i$ is required for all updating schemes; see step \ref{comm_algo} of Alg. \ref{algo_sync}), while first-order methods have lower computational costs. \hide{However, first-order methods suffer from slower convergence speed.} However, a significant reduction of iteration numbers for prescribed accuracy can be achieved by using (quasi) Newton methods. \hide{Moreover, (quasi) Newton methods become necessary when high accuracy solution is desired as first-order methods tend to stagnate at error levels $10^{-2}\sim 10^{-4}$.} \hide{We note that using BFGS estimated curvature information requires less computational resources compared to using the Newton method but the convergence of DRUID-BFGS is inferior to that of the DRUID-Newton, as expected. Therefore, DRUID-BFGS strikes a fine balance between the proposed algorithms with first-order and second-order updates, in terms of computational costs, convergence speed, and achieved solution accuracy.}

\subsection{Distributed Logistic Regression} 
\begin{figure}[t]
  \centering
  \begin{subfigure}[b]{0.49\linewidth}
    \includegraphics[width=\textwidth]{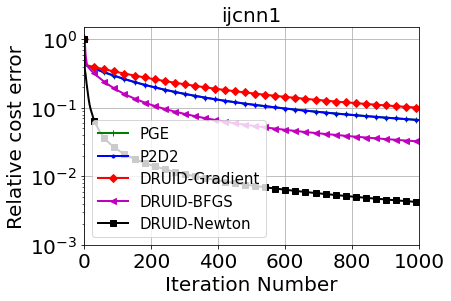}
  \end{subfigure}
  \begin{subfigure}[b]{0.49\linewidth}
    \includegraphics[width=\textwidth]{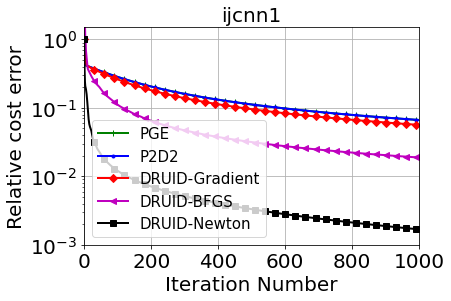}
  \end{subfigure}
  \caption{Performance comparison using the ijcnn1 dataset with different network sizes. We plot the iteration number versus the relative cost error on random graphs with $m=10$ (left) and $m=20$ (right).  }
  \label{logi}
 \end{figure}

The distributed logistic regression solves (\ref{total_obj_function}) with $g(x)=\gamma\norm{x}_1$ and $f_i(\cdot):\mathbb{R}^d\to \mathbb{R}$ defined as:
\begin{gather*}
        f_i(x_i) =  \sum_{j=1}^{m_i} \left[\ln(1+e^{-w_j^\top x_i})+(1-y_j)w_j^\top x_i\right],
\end{gather*}
where $m_i$ is the number of data points accessible by the $i$-th agent. We denote the local training data set as $\{w_j,y_j\}^{m_i}_{j=1}\subset \mathbb{R}^{d}\times \{0,1\}$, where $w_j$ are feature vectors and $y_j$ are known labels. We consider $5,000$ data points from the skin\_noskin dataset with dimensions $d=3$, and $2,000$ data points from the ijcnn1 dataset with dimensions $d=22$. In Figure \ref{sync_async}, we observe that the convergence is slower when only a subset of agents become active (due to less total computation/communication per round compared to the synchronous case). Note that P2D2 and PGE do not support asynchronous implementations. We further explore the effect of the graph topology by varying the size of the network $m$ in Figure \ref{logi}.\hide{As shown in Figure \ref{logi}, algorithms with curvature-aided information consistently outperform first-order algorithms. Moreover,} We observe that DRUID is insensitive to networks with different sizes, but with fixed $p=0.2$. This is consistent with our analysis where the convergence rate is affected by algebraic connectivity, but not system size $m$. 
 
 \begin{figure}[t]
  \centering
  \begin{subfigure}[b]{0.49\linewidth}
    \includegraphics[width=\textwidth]{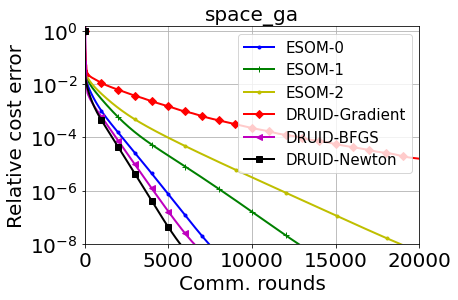}
  \end{subfigure}
  \begin{subfigure}[b]{0.49\linewidth}
    \includegraphics[width=\textwidth]{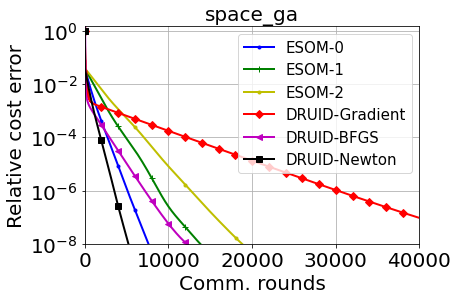}
  \end{subfigure}
  \caption{Performance comparison using the space\_ga dataset. We use the communication rounds as the metric, with the number of agents and probability of generating an edge equal to $m=20,p=0.2$ (left) and $m=40,p=0.8$ (right).}
  \label{ridge}
 \end{figure}
\subsection{Distributed Ridge Regression}
Since existing second-order methods only support differentiable regularization functions, we consider the problem of distributed ridge regression, whose $f_i(\cdot)$ is the same as in (\ref{lasso_f}) but with $g(x)=\gamma \norm{x}^2$. We compare DRUID with ESOM-$K$, where $K$ denotes the number of inner communication loops.\hide{The total number of communication rounds between each agent is used as the metric for comparing, and we note that for ESOM-$0$ (no inner communication loops), the communication costs per round are the same as DRUID.} In the case of ESOM-$K$, a more accurate Hessian estimation can be obtained by increasing $K$, at the cost of more communication rounds. On the other hand, we emphasize that through the use of consensus variables $\{z_{ij}\}$, DRUID-Newton utilizes the exact Hessian without inducing inner loops, and thus achieves the highest communication efficiency. \hide{as shown in the Figure \ref{ridge}.}

\section{Conclusions}

We have proposed a family of distributed primal-dual algorithms for solving convex composite optimization problems. Various computing choices, including gradient, Newton, and BFGS updates, are proposed to achieve a balance between economical computational costs, solution accuracy, and convergence speeds. By use of intermediate consensus variables, we achieve a block-diagonal Hessian that allows us to harness the curvature information without additional communication rounds after each gradient evaluation.\hide{Thus, all methods feature the same minimal communication complexity.} An asynchronous extension of the proposed algorithms is also presented.\hide{that further removes the need for global synchronization.} We establish a unified analytical framework for the proposed algorithms that reveals the difference between various updating schemes. Some future directions include extensions to time-varying and directed network topologies, stochastic gradient evaluation, and hybrid updating schemes.     

\appendices
\section{}\label{AppendixA}
\textit{Proof of Lemma} \ref{lemma1}: The proof was similarly derived in \cite{Mateos2010} and \cite{Wei2014} for the case $g(\cdot)=0$ and the suboptimization problem (\ref{admm1}) was solved exactly. We generalize the results by first writing $y^t=[\alpha^t;\beta^t]$ and recalling the dual update for $y^{t+1}$ in (\ref{admm4}):
$$
    y^{t+1}= y^t+\mu_z(Ax^{t+1}-Bz^{t+1}).
$$ Using (\ref{zupdate}) and premultiplying (\ref{admm4}) with $B^\top$ on both sides, we obtain $B^\top y^{t+1}=0$ for all $t\geq 0$. Since $B^\top = [I_{nd}, I_{nd}]$, it holds that $\alpha^{t+1}+\beta^{t+1}=0$ for all $t\geq 0$. By further assuming $\alpha^0=-\beta^0$, we obtain $\alpha^t=-\beta^t$ for all $t\geq 0$. Recall $A_s:=\hat{A}_s\otimes I_d$ and $A_d:=\hat{A}_d\otimes I_d$, as well as the definition of $A$ in (\ref{prob3}), the dual update (\ref{admm4}) can be rewritten as: 
\begin{align}
    \alpha^{t+1} &=\alpha^t+\mu_z (A_s x^{t+1}-z^{t+1}) ,\label{lemma1_proof_1}\\
    -\alpha^{t+1} & = -\alpha^t+\mu_z(A_dx^{t+1}-z^{t+1}).\label{lemma1_proof_2}
\end{align} 
Recall that $E_s=A_s-A_d$ and $E_u=A_s+A_d$. By taking the sum and difference of (\ref{lemma1_proof_1}) and (\ref{lemma1_proof_2}), we obtain for $t\geq 0$,
\begin{align}
    z^{t+1} &= \tfrac{\mz}{2}E_ux^{t+1}, \label{z_equality}\\
    \alpha^{t+1}&=\alpha^t+\tfrac{\mz}{2}E_sx^{t+1}.\label{lemma1_proof_3}
\end{align}
This establishes that the dual update (\ref{admm4}) for $y^{t+1}$ can be replaced by (\ref{lemma1_proof_3}). Using the definition of $\phi^t=E_s^\top \alpha^t$ and premultiplying (\ref{lemma1_proof_3}) with $E_s^\top$, we obtain (\ref{lemma1_3}). By initializing $z^0=\tfrac{1}{2}E_ux^0$, we have that $z^t=\tfrac{1}{2}E_ux^t$ for $t\geq 0$. Therefore, the update (\ref{zupdate}) for $z^{t+1}$ is not necessary since $z^t$ can be obtained by computing $\tfrac{1}{2}E_ux^t$. It remains to show the equivalence between (\ref{primal}) and (\ref{lemma1_1}). Using (\ref{AL}), it follows that update (\ref{primal}) is given by: 
\begin{gather}
    x^{t+1}=x^t-(H^t)^{-1}\big[\nabla F(x^t)+A^\top y^t+S \lambda^t\nonumber\\
    +\mz A^\top (Ax^t-Bz^t)+\mt S(S^\top x^t-\theta^t)\big] \tag{12}.
\end{gather} Since $y^t= [\alpha^t;-\alpha^t]$ and $z^t=\tfrac{1}{2}E_ux^t$, we obtain:
\begin{align}
    A^\top y^t=\begin{bmatrix}
    A_s^\top & A_d^\top 
    \end{bmatrix}y^t=E_s^\top \alpha^t=\phi^t, \label{lemma1_proof_4}\\
    \mu_z A^\top (Ax^t-Bz^t) = \tfrac{\mu_z}{2}(2D-L_u)x^t=\tfrac{\mu_z}{2}L_sx^t,\label{lemma1_proof_5}
\end{align}
where we have used the identity $A_s^\top -A_d^\top= E_s^\top$, $D=A^\top A$, and $\mz A^\top B z^t=\tfrac{\mz}{2}E_u^\top E_ux^t=\tfrac{\mz}{2}L_u x^t$. After substituting (\ref{lemma1_proof_4}) and (\ref{lemma1_proof_5}) into (\ref{primal}), we obtain the desired.  \QEDB

\textit{Proof of Lemma} \ref{optimal}: The KKT conditions for (\ref{prob3}) are: 
\begin{subequations}
\begin{align}
    \nabla F(x^\star)+A^\top y^\star+S \lambda^\star &= 0,\label{KKT1}\\
    B^\top y^\star &=0,\label{KKT2}\\
    \partial g(\theta^\star)-\lambda^\star &\ni 0, \label{KKT3}\\
    Ax^\star &= Bz^\star ,\label{KKT4}\\
    S^\top x^\star &=\theta^\star.\label{KKT5}
\end{align}
\end{subequations}
Since the objective function is convex with linear constraints, strong duality holds. Recall the definition $B=\begin{bmatrix}
I_{md};I_{md}
\end{bmatrix}\in \mathbb{R}^{2md\times md}$. The condition (\ref{KKT2}) implies that for any dual optimal $y^\star=[\alpha^\star;\beta^\star]$, it holds that $\alpha^\star=-\beta^\star$. Since $A=[A_s;A_d]$ and $E_s=A_s-A_d$, the condition (\ref{KKT1}) can be rewritten as:
\begin{gather}
    \nabla F(x^\star)+ E_s^\top \alpha^\star+S \lambda^\star=0. \label{KKT1_var}
\end{gather}
Note that since $E_s^\top$ has a nontrivial kernel for any network with agent number $m>1$, there exist multiple $\alpha^\star$ that satisfy (\ref{KKT1_var}). We proceed to show that there exists a unique dual optimal $[\alpha^\star;\lambda^\star]$ that lies in the column space of $C=\begin{bmatrix}E_s \\ S^\top \end{bmatrix}$. To show existence, let $\xi^0:=\begin{bmatrix}
\alpha^0;\lambda^0
\end{bmatrix}$ be any dual optimal that satisfies (\ref{KKT1_var}) and (\ref{KKT3}). We denote its projection to the column space of $C$ as $\xi^\star:=[\alpha^\star;\lambda^\star]$. By the property that $C^\top (\xi^0-\xi^\star)=0$, we conclude that $\nabla F(x^\star)+C^\top \xi^\star=0$. Moreover, since $\text{col}(E_s^\top)\,\cap\,\text{col}(S)=0$ and ker$(S)=0$, it holds that $\lambda^0=\lambda^\star$. We prove the uniqueness of $\xi^\star$ by contradiction. Suppose there exist $\xi^1=Cr^1$ and $\xi^2=Cr^2$, $r^1\neq r^2$, that satisfy:
\begin{gather*}
    \nabla F(x^\star)+C^\top Cr^1=0, \\
    \nabla F(x^\star)+C^\top Cr^2=0 .
\end{gather*}
After taking the difference of the above, we obtain $C^\top C(r^1-r^2)=0$. Note that $C^\top C=E_s^\top E_s+S S^\top=L_s+S S^\top$. Since both $L_s$ and $S S^\top$ are positive semidefinite, $C^\top C (r^1-r^2)=0$ if and only if $L_s(r^1-r^2)=S S^\top (r^1-r^2)=0$. Moreover, since the graph is connected, the kernel of $L_s$ is a one dimensional subspace spanned by consensus vector $\mathbf{1}$ and the kernel of $S S^\top$ is spanned by vectors with the $l$-th entry being $0$. Therefore, $L_s(r^1-r^2)=S S^\top(r^1-r^2)=0$ if and only if $r^1-r^2=0$, which contradicts with the assumption $r^1\neq r^2$. \QEDB

\textit{Proof of Lemma} \ref{primal_sub}: Recall the primal update (\ref{lemma1_1}) and the identity $\phi^t=E_s^\top \alpha^t$. After rearranging, we obtain:
\begin{align}
    &\nabla F(x^t)+E_s^\top \alpha^t+S\lambda^t+\tfrac{\mu_z}{2}L_s x^t+\mt S (S^\top x^t-\theta^t)\nonumber\\
    +&H^t(x^{t+1}-x^t)= 0 \label{primal rearrange}
\end{align}
From the dual update (\ref{lemma1_3}), we obtain:
\begin{gather}
    E_s^\top \alpha^t+\tfrac{\mu_z}{2}L_sx^t = E_s^\top \alpha^{t+1}-\tfrac{\mu_z}{2}L_s(x^{t+1}-x^t). \label{E_s alpha}
\end{gather}
Similarly, from the dual update (\ref{lemma1_4}), it holds that
\begin{align}
    &S\lambda^{t}+\mt S(S^\top x^t-\theta^t) \nonumber\\ =&S\lambda^{t+1}-\mt S\big(S^\top(x^{t+1}-x^t)-(\theta^{t+1}-\theta^t)\big). \label{lambda plus mu}
\end{align}
After substituting (\ref{E_s alpha}) and (\ref{lambda plus mu}) into (\ref{primal rearrange}), we obtain:
\begin{align}
    &\nabla F(x^t)+E_s^\top \alpha^{t+1}-\tfrac{\mu_z}{2}L_s(x^{t+1}-x^t)+S\big(\lambda^{t+1}\label{lemma2 proof 28}\\
    -&\mt S^\top(x^{t+1}-x^t)
    +\mt (\theta^{t+1}-\theta^t)\big)+H^t(x^{t+1}-x^t) = 0 \nonumber
\end{align}
Recall $2D=L_s+L_u$ from (\ref{4c}). After adding and subtracting $(\mu_z D+\mt S S^\top+\epsilon I)(x^{t+1}-x^t)$ from (\ref{lemma2 proof 28}), we obtain:
\begin{gather*}
    \nabla F(x^t)+E_s^\top \alpha^{t+1}+\tfrac{\mu_z}{2}L_u(x^{t+1}-x^t)+S\lambda^{t+1}\\
    +\mu_\theta S(\theta^{t+1}-\theta^t)
    +\epsilon(x^{t+1}-x^t)\\
    +\left(H^t-\mu_z D-\mt S S^\top-\epsilon I\right)(x^{t+1}-x^t)=0.
\end{gather*}
Moreover, $\tfrac{\mz}{2}L_u(x^{t+1}-x^t)=\tfrac{\mz}{2}E_u^\top E_u(x^{t+1}-x^t)=\mz E_u^\top (z^{t+1}-z^t)$ by (\ref{z_equality}). Recall $H^t$ in (\ref{H_def}), as well as (\ref{J_def_bfgs}) for the BFGS case. After subtracting KKTa and substituting the definition of $e^t$ in (\ref{error}) and the expression for $\tfrac{\mz}{2}L_u(x^{t+1}-x^t)$ into the above, we obtain the desired. \QEDB 

\section{}\label{AppendixB}

\textit{Proof of Theorem} \ref{sublinear_theorem_combined}: We begin with proving the following two technical inequalities: \begin{align}
       (\lambda^{t+1}-\lambda^t)^\top (\theta^{t+1}-\theta^t) &\geq 0 ,\label{claim1_2}\\
       (\lambda^{t+1}-\lambda^\star)^\top (\theta^{t+1}-\theta^\star) &\geq 0.\label{claim1_1}
\end{align}
From the definition of the proximal operator, it holds that:
\begin{gather}
    \theta^{t+1}=\underset{\theta}{\text{argmin}}\left\{g(\theta)+\tfrac{\mu_\theta}{2}\norm{S^\top x^{t+1}+\tfrac{1}{\mu_\theta}\lambda^t-\theta}^2\right\}. \label{claim1}
\end{gather}
By the optimality condition of (\ref{claim1}), we obtain:
$$
    0\in \partial g(\theta^{t+1})-\mu_\theta(\tfrac{1}{\mu_\theta}\lambda^t+S^\top x^{t+1}-\theta^{t+1})
    =\partial g(\theta^{t+1})-\lambda^{t+1},
$$
where the last equality follows from the dual update (\ref{lemma1_4}). Therefore, it holds that:
$$
    (\lambda^{t+1}-\lambda^t)^\top (\theta^{t+1}-\theta^t)
    \in (\partial g(\theta^{t+1})-\partial g(\theta^t))^\top (\theta^{t+1}-\theta^t)\geq 0 ,
$$
where the inequality follows from the convexity of $g(\cdot)$. Moreover, 
$$
    (\lambda^{t+1}-\lambda^\star)^\top (\theta^{t+1}-\theta^\star)
    \in (\partial g(\theta^{t+1})-\partial g(\theta^\star))^\top (\theta^{t+1}-\theta^\star)\geq 0,
$$
where the inclusion follows from \ref{KKTb}. The rest of the proof is constituted by the following: 
\begin{enumerate}[label=(\roman*)]
    \item Establishing the convergence of $\norm{\va^t-\vas}_{\mathcal{G}^t}^2$ to $0$. 
    \item Establishing the running average upper bound in (\ref{sublinear_theorem_2}). 
\end{enumerate}
Part (i): Since $F(\cdot)$ is convex with the gradient being Lipschitz continuous with parameter $M_f$, the following holds:
\begin{align}
    &\tfrac{1}{M_f}\norm{\nabla F(x^t)-\nabla F(x^\star)}^2 \leq (x^t-x^\star)^\top (\nabla F(x^t)-\nabla F(x^\star)) \nonumber\\
    &=(x^{t+1}-x^\star)^\top (\nabla F(x^t)-\nabla F(x^\star))\nonumber\\
    &\quad+(x^t-x^{t+1})^\top (\nabla F(x^t)-\nabla F(x^\star)).\label{sublinear1}
\end{align}
We proceed to establish an upper bound for the right-hand side of (\ref{sublinear1}) by separately bounding the two components. Recall $H^t$ in (\ref{H_def}). From Lemma \ref{primal_sub}, the following holds:
\begin{align}
    &\nabla F(x^t) - \nabla F(x^\star) = -\big\{E_s^\top( \alpha^{t+1}-\alpha^\star)+ S \Big(\lambda^{t+1}-\lambda^\star\nonumber\\
    &+\mt (\theta^{t+1}-\theta^t)\Big)+\left(J^t+\epsilon I\right)(x^{t+1}-x^t)\nonumber\\
    &+\mz E_u^\top(z^{t+1}-z^t) \big\}\label{F_diff}.
\end{align}
Denoting $\overline{J}^t=J^t+\epsilon I$ and using (\ref{F_diff}), we rewrite the first component of the right hand side of (\ref{sublinear1}) as:
\begin{align}
        &(x^{t+1}-x^\star)^\top (\nabla F(x^t)-\nabla F(x^\star)) 
        = \nonumber\\
    &-(x^{t+1}-x^\star)^\top E_s^\top (\alpha^{t+1}-\alpha^\star)\nonumber\\
    &-(x^{t+1}-x^\star)^\top S \big(\lambda^{t+1}-\lambda^\star+\mt(\theta^{t+1}-\theta^t)\big)\nonumber\\ &-(x^{t+1}-x^\star)^\top \overline{J}^t(x^{t+1}-x^t)\nonumber\\
    &-\mz(x^{t+1}-x^\star)^\top E_u^\top (z^{t+1}-z^t)
   .\label{sublinear2}
\end{align}
From the dual update, KKT conditions, and Lemma \ref{lemma1}, the following holds: 
\begin{align*}
    (x^{t+1}-x^\star)^\top E_s^\top &= \tfrac{2}{\mz}(\alpha^{t+1}-\alpha^t)^\top ,\\
    (x^{t+1}-x^\star)^\top S &= (\theta^{t+1}-\theta^\star)^\top +\tfrac{1}{\mt}(\lambda^{t+1}-\lambda^t)^\top,\\
    (x^{t+1}-x^\star)^\top E_u^\top &= (z^{t+1}-z^\star)^\top.
\end{align*}
Using these expressions for $(x^{t+1}-x^\star)^\top E_s^\top$,  $(x^{t+1}-x^\star)^\top S$, and $(x^{t+1}-x^\star)^\top E_u^\top$, we rewrite (\ref{sublinear2}) as:
\begin{align}
    &(x^{t+1}-x^\star)^\top (\nabla F(x^t)-\nabla F(x^\star)) \label{sublinear6}\\
    &= -\tfrac{2}{\mz}(\alpha^{t+1}-\alpha^t)^\top (\alpha^{t+1}-\alpha^\star)\nonumber\\
    &\underbrace{-(\theta^{t+1}-\theta^\star)^\top (\lambda^{t+1}-\lambda^\star)}_{\leq 0\,\text{from}\,\,(\ref{claim1_1})}-\tfrac{1}{\mt}(\lambda^{t+1}-\lambda^t)^\top(\lambda^{t+1}-\lambda^\star) \nonumber\\
     &-\mt(\theta^{t+1}-\theta^\star)^\top (\theta^{t+1}-\theta^t)\underbrace{-(\lambda^{t+1}-\lambda^t)^\top (\theta^{t+1}-\theta^t)}_{\leq 0\,\,\text{from}\,\,(\ref{claim1_2})} \nonumber\\
    &-(x^{t+1}-x^\star)^\top \overline{J}^t(x^{t+1}-x^t)-2\mz(z^{t+1}-z^\star)^\top (z^{t+1}-z^t)\nonumber\\
    &\leq -\tfrac{2}{\mz}(\alpha^{t+1}-\alpha^t)^\top (\alpha^{t+1}-\alpha^\star)\nonumber\\
    &-\tfrac{1}{\mt}(\lambda^{t+1}-\lambda^t)^\top (\lambda^{t+1}-\lambda^\star)-\mt(\theta^{t+1}-\theta^\star)^\top (\theta^{t+1}-\theta^t)\nonumber\\
    &-(x^{t+1}-x^\star)^\top \overline{J}^t(x^{t+1}-x^t)-2\mz(z^{t+1}-z^\star)^\top(z^{t+1}-z^t)\nonumber\\
    &\underset{\text{(i)}}{=} \tfrac{1}{2}\Big\{\tfrac{2}{\mz}\left(\norm{\alpha^{t}-\alpha^\star}^2-\norm{\alpha^{t+1}-\alpha^\star}^2-\norm{\alpha^{t+1}-\alpha^t}^2\right) \nonumber\\
    &+\tfrac{1}{\mt}\left(\norm{\lambda^{t}-\lambda^\star}^2-\norm{\lambda^{t+1}-\lambda^\star}^2-\norm{\lambda^{t+1}-\lambda^t}^2\right)\nonumber\\
    & +\mt\left(\norm{\theta^{t}-\theta^\star}^2-\norm{\theta^{t+1}-\theta^\star}^2-\norm{\theta^{t+1}-\theta^t}^2\right)\nonumber\\
    &+\norm{x^{t}-x^\star}^2_{\overline{J}^t}-\norm{x^{t+1}-x^\star}^2_{\overline{J}^t}-\norm{x^{t+1}-x^t}^2_{\overline{J}^t}\nonumber\\
    &+2\mz\left(\norm{z^t-z^\star}^2-\norm{z^{t+1}-z^\star}^2-\norm{z^{t+1}-z^t}^2\right)\Big\}\nonumber\\
    &\underset{\text{(ii)}}{=}\tfrac{1}{2}\left(\norm{\va^t-\vas}_{\mathcal{G}^t}^2-\norm{\va^{t+1}-\vas}_{\mathcal{G}^t}^2-\norm{\va^{t+1}-\va^t}_{\mathcal{G}^t}^2\right),\nonumber
\end{align}
where (i) follows from the identity $-2(a-b)^\top  (a-c)=\norm{b-c}^2-\norm{a-b}^2-\norm{a-c}^2$; (ii) follows from the definition (\ref{u_g_def}). We proceed to establish an upper bound for the second term of (\ref{sublinear1}) in the following. Note that $a^\top b \leq \tfrac{1}{2\zeta}a^2+\tfrac{\zeta}{2}b^2$ holds for any $\zeta>0$. By setting $\zeta=\tfrac{M_f}{2},$ we obtain:
\begin{align}
        &(x^t-x^{t+1})^\top (\nabla F(x^t)-\nabla F(x^\star)) \nonumber\\
        &\leq \tfrac{1}{M_f} \norm{\nabla F(x^t)-\nabla F(x^\star)}^2+\tfrac{M_f}{4}\norm{x^{t+1}-x^t}^2.\label{sublinear7}
\end{align}
After substituting (\ref{sublinear6}) and (\ref{sublinear7}) into (\ref{sublinear1}), we obtain: 
\begin{gather*}
        \tfrac{1}{M_f} \norm{\nabla F(x^t)-\nabla F(x^\star)}^2
        \leq 
    \tfrac{1}{2}\Big(\norm{\va^t-\vas}_{\mathcal{G}^t}^2-\norm{\va^{t+1}-\vas}_{\mathcal{G}^t}^2\\
    -\norm{\va^{t+1}-\va^t}_{\mathcal{G}^t}^2 \Big)
   + \tfrac{1}{M_f} \norm{\nabla F(x^t)-\nabla F(x^\star)}^2\\
   +\tfrac{M_f}{4}\norm{x^{t+1}-x^t}^2.
\end{gather*}
By canceling the identical term and rearranging, we obtain: 
\begin{align}
    &\tfrac{1}{2}\left(\norm{\va^t-\vas}^2_{\mathcal{G}^t}-\norm{\va^{t+1}-\vas}^2_{\mathcal{G}^t}\right)\label{sublinear10}\\
    &\geq \tfrac{1}{2}\norm{\va^{t+1}-\va^t}^2_{\mathcal{G}^t}-\tfrac{M_f}{4}\norm{x^{t+1}-x^t}^2 \nonumber\\
    &=\tfrac{1}{2}\big(\norm{x^{t+1}-x^t}^2_{\overline{J}^t-\tfrac{M_f}{2}I}+2\mz\norm{z^{t+1}-z^t}^2 \nonumber\\
    &+\mt\norm{\theta^{t+1}-\theta^t}^2
    +\tfrac{2}{\mz}\norm{\alpha^{t+1}-\alpha^t}^2
    +\tfrac{1}{\mt}\norm{\lambda^{t+1}-\lambda^t}^2\big)\nonumber.
\end{align}
Recall $\overline{J}^t = J^t+\epsilon I$ and we proceed to find a uniform lower bound for $\norm{x^{t+1}-x^t}^2_{\overline{J}^t}$, for gradient descent, Newton, and BFGS computing scheme. Since $J^t_{\mathrm{Gradient}}=0$ and $J^t_{\mathrm{Newton}}=\nabla^2 F(x^t)\succeq0$ by construction, it holds that $\epsilon I\preceq \overline{J}^t_{\mathrm{Gradient/Newton}}$. It remains to find a lower bound for the case of BFGS. Recall $J^t_{\mathrm{BFGS}}$ defined in (\ref{J_def_bfgs}). By the secant condition, $(J^t_{\mathrm{BFGS}}+\mz D+\mt SS^\top +\epsilon I)s^{t-1}=q^{t-1}$, where 
\begin{align*}
    s^{t-1}&=x^{t}-x^{t-1}, \\
    q^{t-1} &= \nabla F(x^{t})-\nabla F(x^{t-1})+\left(\mz D+\mt SS^\top+\epsilon I\right)s^{t-1}.
\end{align*}
Therefore, it holds that: 
\begin{gather}
    J^t_{\mathrm{BFGS}} s^{t-1} =\nabla F(x^{t})-\nabla F(x^{t-1}). \label{modified secant}
\end{gather} 
By premultiplying $(s^{t-1})^\top$ on both sides of (\ref{modified secant}), we obtain: 
$$(s^{t-1})^\top J^t_{\mathrm{BFGS}} s^{t-1} = (x^{t}-x^{t-1})^\top (\nabla F(x^{t})-\nabla F(x^{t-1}))\geq 0.$$
Therefore, the following holds:
$$\sigminwJ\norm{x^{t+1}-x^t}^2\leq \norm{x^{t+1}-x^t}^2_{\overline{J}^t},$$
where $\sigminwJ=\epsilon$. By selecting $\epsilon>\tfrac{M_f}{2}$, we obtain:
$$
    \norm{x^{t+1}-x^t}^2_{\overline{J}^t-\tfrac{M_f}{2}I}\geq \tfrac{\sigminwJ-\tfrac{M_f}{2}}{\sigminwJ}\norm{x^{t+1}-x^t}^2_{\overline{J}^t}.
$$
We denote $\delta= \tfrac{\sigminwJ-\tfrac{M_f}{2}}{\sigminwJ}$ and since $\delta<1$, it holds that:
\begin{gather*}
    2\mz\norm{z^{t+1}-z^t}^2+\mt\norm{\theta^{t+1}-\theta^t}^2+\tfrac{2}{\mz}\norm{\alpha^{t+1}-\alpha^t}^2\\
    +\tfrac{1}{\mt}\norm{\lambda^{t+1}-\lambda^t}^2\geq
    \delta \Big(2\mz\norm{z^{t+1}-z^t}^2+\mt\norm{\theta^{t+1}-\theta^t}^2\\
    +\tfrac{2}{\mz}\norm{\alpha^{t+1}-\alpha^t}^2
    +\tfrac{1}{\mt}\norm{\lambda^{t+1}-\lambda^t}^2\Big).
\end{gather*}
Therefore, (\ref{sublinear10}) can be rewritten as:
    \begin{align}
    \norm{\va^{t+1}-\vas}_{\mathcal{G}^t}^2 \leq \norm{\va^t-\vas}^2_{\mathcal{G}^t}-\delta \norm{\va^{t+1}-\va^t}^2_{\mathcal{G}^t}.\label{sublinear8} 
    \end{align}
Since (\ref{sublinear8}) shows that $\norm{\va^{t}-\vas}^2_{\mathcal{G}^t}$ is monotonically decreasing, it is therefore convergent. We proceed to show Part (ii).

\noindent Part (ii): Recall (\ref{primal rearrange}) and after rearranging, we obtain:
\begin{align}
      &H^t(x^{t+1}-x^t)=-\big\{ \nabla F(x^t)+E_s^\top \alpha^t+S \lambda^t+\tfrac{\mz}{2}L_sx^t \nonumber\\
      &+\mt S (S^\top x^t-\theta^t)\big\}.\label{H_equality}
      \end{align}
Since $H^t=J^t+\mz D+\epsilon I+\mt SS^\top$ as in (\ref{H_def}), an upper bound $H^t\preceq \overline{M}I$ can be obtained by using (\ref{f assumption}) and (\ref{bfgs_bound}):
\begin{align*}
\overline{M}_{\mathrm{Gradient}}&= \mz d_{\mathrm{max}}+\epsilon+\mt,\\
    \overline{M}_{\mathrm{Newton}}&= M_f+\mz d_{\mathrm{max}}+\epsilon+\mt,\\
    \overline{M}_{\mathrm{BFGS}} &= \psi,
\end{align*}
where $d_{\mathrm{max}}=\max_{i} \abs{\mathcal{N}_i}$ denotes the maximum degree. Therefore, the following holds:
\begin{align}
    \overline{M}^2\norm{x^{t+1}-x^t}^2\geq \norm{x^{t+1}-x^t}^2_{(H^t)^2}. \label{sublinear9}
\end{align}
We proceed to establish a lower bound for $\norm{x^{t+1}-x^t}^2_{\overline{J}^t}$:
\begin{align}
    & \norm{x^{t+1}-x^t}^2_{\overline{J}^t}  
    \geq \sigminwJ \norm{x^{t+1}-x^t}^2
    \underset{\text{(i)}}{\geq} \tfrac{\sigminwJ}{\overline{M}^2} \norm{x^{t+1}-x^t}^2_{(H^t)^2} \nonumber\\
    &\underset{\text{(ii)}}{=}\tfrac{\sigminwJ}{\overline{M}^2}\Big\vert\Big\vert\nabla F(x^t)+E_s^\top \alpha^t+S \lambda^t+\tfrac{\mz}{2}L_s x^t
    +\mt S( x^t_l-\theta^t)\Big\vert\Big\vert^2 \nonumber\\
    &\underset{\text{(iii)}}{\geq}\tfrac{\sigminwJ}{\overline{M}^2}\Big( \tfrac{1}{\rho}\norm{\nabla F(x^t)+E_s^\top \alpha^t+S \lambda^t}^2\nonumber\\
    &-\tfrac{1}{\rho-1}\norm{\tfrac{\mz}{2}L_s x^t+\mt S ( S^\top x^t-\theta^t)}^2\Big) \nonumber\\
    &\underset{\text{(iv)}}{\geq}\tfrac{\sigminwJ}{\overline{M}^2}\Big( \tfrac{1}{\rho}\norm{\nabla F(x^t)+E_s^\top \alpha^t+S \lambda^t}^2 - \tfrac{2}{\rho-1}\norm{\tfrac{\mz}{2}L_sx^t}^2\nonumber\\
    &\qquad-\tfrac{2}{\rho-1}\norm{\mt  (x^t_l-\theta^t)}^2\Big),\label{sublinear_x}
\end{align}
where (i) follows from (\ref{sublinear9}); (ii) follows from (\ref{H_equality}); (iii) follows from $(a+b)^2\geq \tfrac{1}{\rho}a^2-\tfrac{1}{\rho-1}b^2$ for any $\rho>1$; (iv) follows from $-(a+b)^2\geq -2(a^2+b^2)$. Also note that $\norm{\mt S(S^\top x^t-\theta^t)}=\norm{\mt(x^t_l-\theta^t)}$ by definition of $S=s_l\otimes I_d$ being the selection matrix. Further observe that the following holds due to dual updates (\ref{lemma1_3}) and (\ref{lemma1_4}):
\begin{align*}
        \alpha^{t+1}-\alpha^t &= \tfrac{\mz}{2} E_s x^{t+1},\\
    \lambda^{t+1}-\lambda^t   &=\mt (S^\top x^{t+1}-\theta^{t+1}).
\end{align*}
Therefore, we obtain the following:
\begin{align}
        \tfrac{2}{\mz}\norm{\alpha^{t+1}-\alpha^t}^2 &= \tfrac{\mz}{2}\norm{E_sx^{t+1}}^2=\tfrac{\mz}{2}\norm{x^{t+1}}^2_{L_s} ,\label{dual_1}\\
    \tfrac{1}{\mt}\norm{\lambda^{t+1}-\lambda^t}^2&=\mt\norm{x^{t+1}_l-\theta^{t+1}}^2,\label{dual_2}
\end{align}
By denoting the maximum eigenvalue of $L_s$ as $\sigmaxLs$ and selecting $\rho-1>\sigmaxLs$, we obtain:
\begin{align}
    \tfrac{\sigminwJ}{\overline{M}^2}\tfrac{2}{\rho-1}\norm{\tfrac{\mz}{2}L_sx^t}^2\leq \tfrac{\sigminwJ \mz^2}{2\overline{M}^2\sigmaxLs}\norm{x^t}^2_{(L_s)^2}\leq \tfrac{\sigminwJ \mz^2}{2\overline{M}^2}\norm{x^t}^2_{L_s}.\label{x_Ls_bd}
\end{align}
Recall the definition (\ref{u_g_def}). We establish (\ref{sublinear_theorem_2}) as follows: 
\begin{align*}
      &\tfrac{1}{T}\tfrac{\mz}{2}\norm{x^1}^2_{L_s}+\tfrac{\mt}{T}\norm{x^1_l-\theta^1}^2+\tfrac{1}{T}\sum_{t=1}^T\norm{\va^{t+1}-\va^t}^2_{\mathcal{G}^t}\\
      &=\tfrac{1}{T}\tfrac{\mz}{2}\norm{x^1}_{L_s}^2+\tfrac{\mt}{T}\norm{x^1_l-\theta^1}^2+\tfrac{1}{T}\sum_{t=1}^T\bigg(\norm{x^{t+1}-x^t}^2_{\overline{J}^t}\\
      &+2\mz\norm{z^{t+1}-z^t}^2+\mt\norm{\theta^{t+1}-\theta^t}^2
      +\tfrac{2}{\mz}\norm{\alpha^{t+1}-\alpha^t}^2\\
      &+\tfrac{1}{\mt}\norm{\lambda^{t+1}-\lambda^t}^2\bigg)
      \underset{\text{(i)}}{\geq}\tfrac{1}{T}\tfrac{\mz}{2}\norm{x^{T+1}}^2_{L_s}+\tfrac{\mt}{T}\norm{x^{T+1}_l-\theta^{T+1}}^2\\
      &+ \tfrac{1}{T}\sum_{t=1}^T\bigg(\tfrac{\sigminwJ}{\overline{M}^2\rho}\norm{\nabla F(x^t)+E_s^\top \alpha^t+S \lambda^t}^2+2\mz\norm{z^{t+1}-z^t}^2\\
      &+\mt\norm{\theta^{t+1}-\theta^t}^2
      +\left(\tfrac{\mz}{2}-\tfrac{\sigminwJ \mz^2}{2\overline{M}^2}\right)\norm{x^t}^2_{L_s}\\
      &+\left(\mt -\tfrac{2\sigminwJ \mt^2}{\overline{M}^2(\rho-1)}\right)\norm{x^t_l-\theta^t}^2\bigg)
\end{align*}
where (i) follows from substituting (\ref{sublinear_x})-(\ref{x_Ls_bd}). All coefficients are ensured to be positive by selecting: $\mz\epsilon<\psi^2$, and $\rho > \max\left\{\tfrac{2\sigminwJ \mt}{\overline{M}^2},\sigmaxLs\right\}+1$, where $\sigma^{\overline{J}}_{\mathrm{min}}=\epsilon$.

\textit{Proof of Corollary} \ref{cor1}: Following Theorem \ref{sublinear_theorem_combined} and standard analysis techniques in \cite{He1994} and \cite{Yin2017}, we obtain that $\norm{\va^{t}-\vas}\to 0$ as $t\to \infty$. After taking telescoping sum from $t=1$ to $\infty$ on both sides of (\ref{sublinear8}), we obtain:
$$
    \delta \sum_{t=1}^\infty \norm{\va^{t+1}-\va^t}^2_{\mathcal{G}^t}\leq \norm{\va^1-\vas}^2_{\mathcal{G}^t}, 
$$
i.e., $\sum_{t=1}^\infty \norm{\va^{t+1}-\va^t}^2_{\mathcal{G}^t}$ is bounded. Define $b^T:=\tfrac{1}{T}\sum_{t=1}^T \norm{\va^{t+1}-\va^t}^2_{\mathcal{G}^t}$. Then $\lim_{T\to \infty} Tb^T=\lim_{T\to\infty}\sum_{t=1}^T\norm{\va^{t+1}-\va^t}^2_{\mathcal{G}^t}<\infty$. Therefore, $b^T=\tfrac{1}{T}\sum_{t=1}^T \norm{\va^{t+1}-\va^t}^2_{\mathcal{G}^t}=\mathcal{O}(\tfrac{1}{T}).$ By (\ref{sublinear_theorem_2}), each term in (\ref{corollary1}) is of order $\mathcal{O}(\tfrac{1}{T})$. \QEDB

\section{}\label{AppendixC}
\textit{Proof of Lemma} \ref{lemma_error_bound}: Recall the definition of $e^t$ in (\ref{error}):
$$
    e^t=\nabla F(x^t)-\nabla F(x^{t+1})+J^t(x^{t+1}-x^t).
$$
By applying the triangle and Cauchy-Schwartz inequality, we obtain:
\begin{align}
    \norm{e^t}\leq \norm{\nabla F(x^t)-\nabla F(x^{t+1})}+\norm{J^t}\norm{x^{t+1}-x^t}.\label{lemma5_1}
\end{align}
In the case of gradient updates, $J^t=0$. Therefore, 
$$
   \norm{ e^t_{\mathrm{Gradient}}}\leq \norm{\nabla F(x^t)-\nabla F(x^{t+1})}
   \leq M_f\norm{x^{t+1}-x^t},
$$
where the last inequality follows from Assumption \ref{assumption1}. Setting $\tau^t_{\mathrm{Gradient}}=M_f$, we obtain (\ref{error_bound_gradient}). In the case of Newton updates, $J^t=\nabla^2 F(x^t)$. By Assumption \ref{assumption1} and (\ref{lemma5_1}), we obtain:
\begin{align}
    \norm{e^t}\leq 2M_f\norm{x^{t+1}-x^t}. \label{lemma5_4}
\end{align}
Moreover, by the fundamental theorem of calculus, $\nabla F(x^{t+1})-\nabla F(x^t)$ can be written as: 
$$
    \nabla F(x^{t+1})-\nabla F(x^t)
    =\int_0^1 \nabla^2 F(sx^{t+1}+(1-s)x^t)(x^{t+1}-x^t) ds. 
$$
By adding and subtracting $\int_0^1 \nabla^2 F(x^t)(x^{t+1}-x^t)ds$, we further obtain:
\begin{gather*}
    \nabla F(x^{t+1})-\nabla F(x^t)
    =\int_0^1 \nabla^2 F(x^t)(x^{t+1}-x^t) ds\\
    +\int_0^1 \left(\nabla^2 F(sx^{t+1}+(1-s)x^t)-\nabla^2 F(x^t)\right)(x^{t+1}-x^t) ds.
\end{gather*}
Since the integrand of the first term is constant with respect to $s$, it holds that:
\begin{align*}
    &\norm{\nabla F(x^{t+1})-\nabla F(x^t)-\nabla^2 F(x^t)(x^{t+1}-x^t)} 
    =\\
    &\norm{\int_0^1 \left(\nabla^2 F(sx^{t+1}+(1-s)x^t)-\nabla^2 F(x^t)\right)(x^{t+1}-x^t)ds}
    \leq\\
    &\int_0^1 \norm{\nabla ^2 F(sx^{t+1}+(1-s)x^t)-\nabla^2 F(x^t)}\cdot\norm{x^{t+1}-x^t} ds 
    \leq\\
    &\int_0^1 sL_f\norm{x^{t+1}-x^t}^2 ds
    =\tfrac{L_f}{2}\norm{x^{t+1}-x^t}^2. 
\end{align*}
Note that in the case of Newton updates, $$\norm{e^t}=\norm{\nabla F(x^{t+1})-\nabla F(x^t)-\nabla^2 F(x^t)(x^{t+1}-x^t)}.$$ By combining (\ref{lemma5_4}) and the above, we obtain:
$
    \norm{e^t}\leq \tau^t_{\mathrm{Newton}}\norm{x^{t+1}-x^t},
$ where $\tau^t_\mathrm{Newton}$ is defined in (\ref{error_bound_newton}). We proceed to establish (\ref{error_bound_bfgs}). Recall the definition of $J^t_{\mathrm{BFGS}}$ in (\ref{J_def_bfgs}):
\begin{gather}
    J^t_{\mathrm{BFGS}}= H^t_{\mathrm{BFGS}}-\mz D-\mt SS^\top-\epsilon I_{md}.\label{lemma5_J}
\end{gather}
Therefore, $H^{t+1}$ (suppressing the subscript BFGS) satisfies the secant condition: 
$
    H^{t+1} s^t = q^t, 
$
where $\{q^t,s^t\}$ as per the definition in (\ref{diff}) can be written as:
\begin{align*}
    s^t &=x^{t+1}-x^t,\\
    q^t 
    &=\nabla F(x^{t+1})-\nabla F(x^t)+\big(\mz D+\mt S S^\top+\epsilon I \big)s^t.
\end{align*}
From the secant condition, it holds that:
\begin{gather*}
    \nabla F(x^t)-\nabla F(x^{t+1})
    =\\
    -\left(H^{t+1}-\mz D-\mt S S^\top-\epsilon I\right) (x^{t+1}-x^t). 
\end{gather*}
Using (\ref{lemma5_J}) and the expression for $\nabla F(x^t)-\nabla F(x^{t+1})$ into (\ref{error}), we obtain:
\begin{align*}
    \norm{e^t} &= \norm{\left(H^t-H^{t+1}\right)(x^{t+1}-x^t)} \\
    &\leq \norm{H^t-H^{t+1}}\norm{x^{t+1}-x^t}. 
\end{align*}
Denoting $\tau^t_{\mathrm{BFGS}}=\norm{H^{t}-H^{t+1}}$ and using (\ref{bfgs_bound}), we obtain (\ref{error_bound_bfgs}). \QEDB 

The following Lemma that will be useful for establishing Theorem \ref{theorem_linear_syn}. 

\begin{lemma}\label{lemma6} Recall $C:=$ $\begin{bmatrix}
E_s \\ S^\top
\end{bmatrix}$ and $\phi^t=E_s^\top \alpha^t$ in (\ref{updates}). Denote the smallest positive eigenvalue of $C C^\top$ as $\sigma^+_{\mathrm{min}}$ and consider the unique dual optimal pair $(\alpha^\star,\lambda^\star)$ that lies in the column space of $C$ as established in Lemma \ref{optimal}. The following holds: 
\begin{align}
    &\sigma^+_{\mathrm{min}}\left(\norm{\alpha^{t+1}-\alpha^\star}^2+\norm{\lambda^{t+1}-\lambda^\star}^2\right)\nonumber\\
    &\leq \norm{E_s^\top (\alpha^{t+1}-\alpha^\star)+S(\lambda^{t+1}-\lambda^\star) }^2.\label{lemma_col_space}
\end{align}
\textit{Proof}: We proceed by showing that $[\alpha^{t+1};\lambda^{t+1}]$ lies in $\mathrm{col}(C)$. We rewrite dual updates (\ref{lemma1_3})--(\ref{lemma1_4})  as:
$$
    \begin{bmatrix}
    \alpha^{t+1}\\
    \lambda^{t+1}
    \end{bmatrix}=
    \begin{bmatrix}
    \alpha^{t}\\
    \lambda^{t}
    \end{bmatrix}+
    \begin{bmatrix}
    \tfrac{\mz}{2}E_s\\
    \mt S^\top
    \end{bmatrix}x^{t+1}-
    \begin{bmatrix}
    0 \\ 
    \mt I_d
    \end{bmatrix}\theta^{t+1}.
$$
We show that the column space of $M:=\begin{bmatrix}
0\\ \mt I_d
\end{bmatrix}$ belongs in the column space of $N:=\begin{bmatrix}
\tfrac{\mz}{2}E_s \\ \mt S^\top
\end{bmatrix}$. Consider fixed $r^x\in\mathbb{R}^d$. Let $r^y\in\mathbb{R}^{md}$ such that each sub-vector component $r^y_i=r^x$, i.e., $r^y=[r^x;\dots;r^x]$. Then it holds that 
$$
    \begin{bmatrix}
\tfrac{\mz}{2}E_s \\ \mt S^\top
\end{bmatrix}r^y = 
\begin{bmatrix}
0\\
\mt r^y_l
\end{bmatrix}=
\begin{bmatrix}
0 \\ \mt I_d
\end{bmatrix}r^x.
$$, which shows $\mathrm{col}(M)\subset \mathrm{col}(N)$. By choosing $\mz=2\mt$, we conclude that $[\alpha^{t+1}-\alpha^\star;\lambda^{t+1}-\lambda^\star]$ lies in the column space of $C$.  \QEDB
\end{lemma}

\textit{Proof of Theorem }\ref{theorem_linear_syn}: Using Lemma \ref{primal_sub}, we obtain:
\begin{gather*}
    \nabla F(x^{t+1})-\nabla F(x^\star) = -\big(E_s^\top (\alpha^{t+1}-\alpha^\star)+\epsilon(x^{t+1}-x^t)\\
    +S(\lambda^{t+1}-\lambda^\star+\mt (\theta^{t+1}-\theta^t))+e^t+\mz E_u^\top(z^{t+1}-z^t)), 
\end{gather*}
Since $F(x)$ is strongly convex with Lipschitz continuous gradient, the following inequality holds \cite{Nesterov2018}:
\begin{gather*}
    \tfrac{m_fM_f}{m_f+M_f}\norm{x^{t+1}-x^\star}^2+\tfrac{1}{m_f+M_f}\norm{\nabla F(x^{t+1})-\nabla F(x^\star)}^2 
    \leq\\
    (x^{t+1}-x^\star)^\top (\nabla F(x^{t+1})-\nabla F(x^\star)) 
    \leq\\
    -(x^{t+1}-x^\star)^\top e^t
    -\epsilon(x^{t+1}-x^\star)^\top(x^{t+1}-x^t) \\
    -(x^{t+1}-x^\star)E_s^\top (\alpha^{t+1}-\alpha^\star)
    -(x^{t+1}-x^\star)^\top S \Big(\lambda^{t+1}-\lambda^\star\\
    +\mt(\theta^{t+1}-\theta^t)\Big)
    -\mz(x^{t+1}-x^\star)^\top E_u^\top (z^{t+1}-z^t), \nonumber\label{main_1}
\end{gather*}
where the last inequality follows from substituting the expression of $\nabla F(x^{t+1})-\nabla F(x^\star)$ above. Using similar techniques used in deriving (\ref{sublinear1})-(\ref{sublinear6}), we obtain
\begin{align}
       &\tfrac{2m_fM_f}{m_f+M_f}\norm{x^{t+1}-x^\star}^2+\tfrac{2}{m_f+M_f}\norm{\nabla F(x^{t+1})-\nabla F(x^\star)}^2\nonumber\\
       &\leq  \epsilon\big(\norm{x^t-x^\star}^2-\norm{x^{t+1}-x^\star}^2-\norm{x^{t+1}-x^t}^2\big)\nonumber\\
       & +2\mz\big(\norm{z^t-z^\star}^2-\norm{z^{t+1}-z^\star}^2-\norm{z^{t+1}-z^t}^2\big)\nonumber\\
       & +\tfrac{1}{\mu_\theta}\big(\norm{\lambda^t-\lambda^\star}^2-\norm{\lambda^{t+1}-\lambda^\star}^2-\norm{\lambda^{t+1}-\lambda^t}^2\big)\nonumber\\
       & +\mt\big(\norm{\theta^t-\theta^\star}^2-\norm{\theta^{t+1}-\theta^\star}^2-\norm{\theta^{t+1}-\theta^t}^2\big)\nonumber\\
       & +\tfrac{1}{\mz}\big(\norm{\alpha^{t+1}-\alpha^\star}-\norm{\alpha^{t+1}-\alpha^\star}-\norm{\alpha^{t+1}-\alpha^t}^2\big)\nonumber\\
       & -2(x^{t+1}-x^\star)^\top e^t \nonumber \\
       &=\norm{\va^t-\vas}^2_{\mathcal{H}}-\norm{\va^{t+1}-\vas}^2_\mathcal{H}-\norm{\va^{t+1}-\va^t}^2_{\mathcal{H}}\nonumber\\
       & -2(x^{t+1}-x^\star)^\top e^t,
\end{align}
Therefore, we obtain:
\begin{align}
    &\tfrac{2m_fM_f}{m_f+M_f}\norm{x^{t+1}-x^\star}^2+\tfrac{2}{m_f+M_f}\norm{\nabla F(x^{t+1})-\nabla F(x^\star)}^2\nonumber \\
    &\quad +\norm{\va^{t+1}-\va^t}^2_{\mathcal{H}}+2(x^{t+1}-x^\star)^\top e^t \nonumber\\
    &\leq \norm{\va^t-\vas}^2_{\mathcal{H}}-\norm{\va^{t+1}-\vas}^2_{\mathcal{H}}.  \label{theorem_lb}
\end{align}
To establish linear convergence, we need to show the following holds for some $\eta>0$:
\begin{gather}
    \eta \norm{\va^{t+1}-\vas}^2_{\mathcal{H}}\leq \norm{\va^t-\vas}^2_{\mathcal{H}}-\norm{\va^{t+1}-\vas}^2_{\mathcal{H}}. \label{suffice2}
\end{gather}
We expand the expression of $\eta\norm{\va^{t+1}-\vas}^2_{\mathcal{H}}$ as follows:
\begin{gather}
    \eta \norm{\va^{t+1}-\vas}^2_{\mathcal{H}}
    = \eta\Big(\epsilon\norm{x^{t+1}-x^\star}^2+2\mz\norm{z^{t+1}-z^\star}^2\nonumber\\
    +\tfrac{2}{\mz}\norm{\alpha^{t+1}-\alpha^\star}^2 
    +\mt\norm{\theta^{t+1}-\theta^\star}^2
    +\tfrac{1}{\mt}\norm{\lambda^{t+1}-\lambda^\star}^2\Big).\label{theorem_linear_component}
\end{gather}
We proceed to establish an upper bound for each component of (\ref{theorem_linear_component}). From Lemma \ref{primal_sub}, the following holds:
\begin{gather*}
     E_s^\top(\alpha^{t+1}-\alpha^\star)+S (\lambda^{t+1}-\lambda^\star)
    =-\Big\{\nabla F(x^{t+1})-\nabla F(x^\star)\\
    +\epsilon(x^{t+1}-x^t)
    +\mz E_u^\top(z^{t+1}-z^t)+\mt S(\theta^{t+1}-\theta^t)+e^t\Big\}.
\end{gather*}
Then we obtain:
\begin{align}
    &\sigma^+_{\mathrm{min}}\left(\norm{\alpha^{t+1}-\alpha^\star}^2+\norm{\lambda^{t+1}-\lambda^\star}^2\right)\nonumber\\
    &\underset{\text{(i)}}{\leq} \norm{E_s^\top (\alpha^{t+1}-\alpha^\star)+S (\lambda^{t+1}-\lambda^\star)}^2\nonumber\\
    &\underset{\text{(ii)}}{\leq} 5\Big(\norm{\nabla F(x^{t+1})-\nabla F(x^\star)}^2+\epsilon^2\norm{x^{t+1}-x^t}^2\nonumber\\
    &+\mt^2\norm{\theta^{t+1}-\theta^t}^2+\norm{e^t}^2+\sigma^{L_u}_{\mathrm{max}}\mz^2\norm{z^{t+1}-z^t}^2\Big), \label{theorem_linear_dual_bound}
\end{align}
where (i) follows from Lemma \ref{lemma6}; (ii) follows from the inequality $(\sum_{i=1}^n a_i)^2\leq \sum_{i=1}^n n a_i^2$. Recalling that we have selected $\mz=2\mt$, we obtain:
\begin{align}
    &\tfrac{2}{\mz}\norm{\alpha^{t+1}-\alpha^\star}^2+\tfrac{1}{\mt} \norm{\lambda^{t+1}-\lambda^\star}^2 \nonumber\\
    &= \tfrac{1}{\mt} \left(\norm{\alpha^{t+1}-\alpha^\star}^2+\norm{\lambda^{t+1}-\lambda^\star}^2\right)\nonumber\\
    &\underset{\text{(i)}}{\leq} \tfrac{5}{\mt\sigma^+_{\mathrm{min}}}\Big(\norm{\nabla F(x^{t+1})-\nabla F(x^\star)}^2+\epsilon^2\norm{x^{t+1}-x^t}^2\nonumber\nonumber\\
    &+\mt^2\norm{\theta^{t+1}-\theta^t}^2+\norm{e^t}^2+\sigma^{L_u}_{\mathrm{max}}\mz^2\norm{z^{t+1}-z^t}^2\Big),\label{theorem_linear_part1}
\end{align}
where (i) follows from dividing (\ref{theorem_linear_dual_bound}) by $\sigma^+_{\mathrm{min}}$ on both sides and substituting. Note that since $z^{t+1}-z^\star= \tfrac{1}{2}E_u(x^{t+1}-x^\star)$, it holds that: $$2\mz\norm{z^{t+1}-z^\star}^2\leq \tfrac{\mz\sigma^{L_u}_{\mathrm{max}}}{2}\norm{x^{t+1}-x^\star}^2.$$
Using the upper bound for $2\mz\norm{z^{t+1}-z^\star}^2$, the inequality (\ref{theorem_linear_part1}), and $\mt\norm{\theta^{t+1}-\theta^\star}^2\leq 2\mt \norm{x^{t+1}-x^\star}^2+\tfrac{2}{\mt} \norm{\lambda^{t+1}-\lambda^t}^2$ from (\ref{lemma1_4}) and KKTd, we obtain an upper bound for (\ref{theorem_linear_component}) as:
\begin{gather*}
    \eta \norm{\va^{t+1}-\vas}^2_{\mathcal{H}}\leq \eta\Big\{\tfrac{5}{\mt\sigma^+_{\mathrm{min}}}\Big(\norm{\nabla F(x^{t+1})-\nabla F(x^\star)}^2\\
    +\epsilon^2\norm{x^{t+1}-x^t}^2
    +\mt^2\norm{\theta^{t+1}-\theta^t}^2+\norm{e^t}^2\\
    +\sigma^{L_u}_{\mathrm{max}}\mz^2\norm{z^{t+1}-z^t}^2\Big)+\tfrac{2}{\mt}\norm{\lambda^{t+1}-\lambda^t}^2\\
    +(\epsilon+2\mt+\tfrac{\mz\sigma^{L_u}_{\mathrm{max}}}{2})\norm{x^{t+1}-x^\star}^2\Big\}. 
\end{gather*}
Recall that the right-hand side of (\ref{suffice2}) is lower bounded as in (\ref{theorem_lb}). Therefore, it suffices to  prove the following to establish (\ref{suffice2}): 
\begin{align}
    &\eta\Big\{\tfrac{5}{\mt\sigma^+_{\mathrm{min}}}\Big(\norm{\nabla F(x^{t+1})-\nabla F(x^\star)}^2+\epsilon^2\norm{x^{t+1}-x^t}^2\nonumber\\
    &+\mt^2\norm{\theta^{t+1}-\theta^t}^2+\norm{e^t}^2+\sigma^{L_u}_{\mathrm{max}}\mz^2\norm{z^{t+1}-z^t}^2\Big)\nonumber\\
    &+\tfrac{2}{\mt}\norm{\lambda^{t+1}-\lambda^t}^2+(\epsilon+2\mt+\tfrac{\mz\sigma^{L_u}_{\mathrm{max}}}{2})\norm{x^{t+1}-x^\star}^2\Big\}\nonumber\\
    &\leq \tfrac{2m_fM_f}{m_f+M_f}\norm{x^{t+1}-x^\star}^2+\tfrac{2}{m_f+M_f}\norm{\nabla F(x^{t+1})-\nabla F(x^\star)}^2\nonumber \\
    &+\norm{\va^{t+1}-\va^t}^2_{\mathcal{H}}+2(x^{t+1}-x^\star)^\top e^t,\label{main_2}
\end{align}
Note that $-\zeta \norm{e^t}^2-\tfrac{1}{\zeta}\norm{x^{t+1}-x^\star}^2\leq 2(x^{t+1}-x^\star)^\top e^t$ holds for any $\zeta>0$. To prove (\ref{main_2}), it is therefore sufficient to show:
\begin{align}
      &\zeta(\tau^t)^2\norm{x^{t+1}-x^t}^2+\eta\Big\{\tfrac{5}{\mt\sigma^+_{\mathrm{min}}}\Big(\norm{\nabla F(x^{t+1})-\nabla F(x^\star)}^2\nonumber\\
      &+((\tau^t)^2+\epsilon^2)\norm{x^{t+1}-x^t}^2
    +\mt^2\norm{\theta^{t+1}-\theta^t}^2\nonumber\\
    &+\sigma^{L_u}_{\mathrm{max}}\mz^2\norm{z^{t+1}-z^t}^2\Big)+\tfrac{2}{\mt}\norm{\lambda^{t+1}-\lambda^t}^2\nonumber\\
    &+(\epsilon+2\mt+\tfrac{\mz\sigma^{L_u}_{\mathrm{max}}}{2})\norm{x^{t+1}-x^\star}^2\Big\}\label{final}\\
    &\leq\left(\tfrac{2m_fM_f}{m_f+M_f}-\tfrac{1}{\zeta}\right)\norm{x^{t+1}-x^\star}^2
    +\epsilon\norm{x^{t+1}-x^t}^2\nonumber\\
    &+2\mz\norm{z^{t+1}-z^t}^2
    +\tfrac{2}{\mz}\norm{\alpha^{t+1}-\alpha^t}^2+\mt\norm{\theta^{t+1}-\theta^t}^2\nonumber\\
    &+\tfrac{1}{\mt}\norm{\lambda^{t+1}-\lambda^t}^2
    +\tfrac{2}{m_f+M_f}\norm{\nabla F(x^{t+1})-\nabla F(x^\star)}^2 \nonumber
\end{align}
where we have used $\norm{e^t}^2\leq (\tau^t)^2 \norm{x^{t+1}-x^t}^2$ from Lemma \ref{lemma_error_bound}. Establishing (\ref{final}) amounts to ensuring the coefficient of each term in the left-hand side is bounded by the coefficient of the corresponding term on the right-hand side. By selecting $\eta$ as in (\ref{eta}), we establish (\ref{final}). Therefore, the inequality (\ref{suffice2}) holds, which equivalently establishes the linear convergence rate. \QEDB

\textit{Proof of Theorem }\ref{theorem_asy}: The proof proceeds as follows:
\begin{align}
    &\norm{\va^{t+1}-\vas}^2_{\mathcal{H}\Omega_\alpha^{-1}}
    = \norm{\va^{t}+\Omega_\alpha^{t+1}(T\va^t-\va^t)-\vas}^2_{\mathcal{H}\Omega_\alpha^{-1}} \nonumber \\
    &= \norm{\va^t-\vas}_{\mathcal{H}\Omega_\alpha^{-1}}^2 +2(\va^t-\vas)^\top\mathcal{H}\Omega_\alpha^{-1}\Omega_\alpha^{t+1}(T\va^t-\va^t)\nonumber\\ & + (T\va^t-\va^t)^\top\Omega_\alpha^{t+1}\mathcal{H}\Omega_\alpha^{-1}\Omega_\alpha^{t+1}(T\va^t-\va^t), \label{theorem_linear_asy}
\end{align}
Since $\Omega_\alpha^{t+1},\Omega^{-1}_\alpha$, and $\mathcal{H}$ are all diagonal matrices, they commute with each other. Moreover, since each sub-block of $\Omega_\alpha^{t+1}$ is $I_d$ or 0, it holds that $\Omega_\alpha^{t+1}\Omega_\alpha^{t+1}=\Omega^{t+1}_\alpha$. After taking conditional expectation on both sides of (\ref{theorem_linear_asy}), we obtain:
\begin{gather*}
     \mathbb{E}^t\left[\norm{\va^{t+1}-\vas}^2_{\mathcal{H}\Omega_\alpha^{-1}} \right] 
     = \norm{\va^{t}-\vas}_{\mathcal{H}\Omega_\alpha^{-1}}^2
+\norm{T\va^t-\va^t}^2_{\mathcal{H}}\\
+2(\va^t-\vas)^\top\mathcal{H}(T\va^t-\va^t)  
\underset{\text{(i)}}{\leq}\\
\norm{\va^t-\vas}_{\mathcal{H}\Omega_\alpha^{-1}}^2-\tfrac{\eta}{1+\eta}\norm{\va^t-\vas}^2_{\mathcal{H}} 
\underset{\text{(ii)}}{\leq}\\
\left(1-\tfrac{p^{\mathrm{min}}\eta }{1+\eta} \right)\norm{\va^t-\vas}^2_{\mathcal{H}\Omega_\alpha^{-1}},
\end{gather*}
where (i) follows from the fact that $ 2(\va-\vas)^\top\mathcal{H}(T\va-\va)+\norm{T\va-\va}_{\mathcal{H}}^2 \leq  -\tfrac{\eta}{1+\eta}\norm{\va-\vas}^2_{\mathcal{H}}$ holds for any $\va\in\mathbb{R}^{(m+2n+2)d}$ using Theorem \ref{theorem_linear_syn}; (ii) follows from $\tfrac{\eta}{1+\eta}\norm{\va^t-\vas}_{\mathcal{H}}\geq \tfrac{p^{\mathrm{min}}\eta}{1+\eta}\norm{\va^t-\vas}^2_{\mathcal{H}\Omega_\alpha^{-1}}$.  \QEDB

\textit{Proof of Corollary} \ref{corollary_asy}: We first distribute each $\alpha_k,k\in[n]$, to each edge and label agents and edges with an arbitrary order. For each edge $\mathcal{E}_k$, we write $\mathcal{E}_k=(i,j)$ with the convention $i<j$. For each agent $i$, we divide the incident edges to two groups:  $\mathcal{P}_i=\{k:\mathcal{E}_k=(i,j),j\in\mathcal{N}_i\}$ and $\mathcal{S}_i=\{k:\mathcal{E}_k=(j,i),j\in\mathcal{N}_i\}$. Consider the activation scheme using $\Omega^{t+1}$. Recall $\alpha^{t+1}_k=\alpha^t_k+\tfrac{\mz}{2}(x_i^{t+1}-x_j^{t+1})$. The dual updates are described by:
\begin{align*}
    \phi^{t+1}_i &= \phi^t_i +\tfrac{\mz}{2}X^{t+1}_{ii}\sum_{j\in\mathcal{N}_i}(x_i^{t+1}-x_j^{t+1})\\
    &= \phi_i^t +X^{t+1}_{ii} \left\{\sum_{k\in \mathcal{P}_i}(\alpha^{t+1}_k-\alpha^{t}_k)+\sum_{k\in \mathcal{S}_i}(\alpha^{t}_k-\alpha^{t+1}_k)\right\}.
\end{align*}
Therefore, if $X^{t+1}_{ii}=I_d$, then $Y^{t+1}_{kk}=I_d$ for $k\in P_i\cup S_i$ for the corresponding $\Omega^{t+1}_\alpha$, i.e., all incident edges are active. It can be verified that we can map $X^{t+1}$ to $Y^{t+1}$ as:
$$
    Y^{t+1} = \textbf{Blkdiag}\left(\ceil{\frac{E_u X^{t+1} (\mathbf{1}\otimes I_d)}{2}} \right),
$$
where $\lceil\cdot \rceil$ is the entry-wise ceiling operation and $\mathbf{1}\in\mathbb{R}^{m}$ is the all one vector. To show $\mathbb{E}^t[\Omega^{t+1}_\alpha]\succ0$, we only need to show $\mathbb{E}^t[Y^{t+1}]\succ 0$, which amounts to showing that $\mathbb{E}^t\left[\ceil{\frac{E_u X^{t+1} (\mathbf{1}\otimes I_d)}{2}}_{k}\right]\in \mathbb{R}^{d\times d},k\in[n],$ is positive definite. Note that:
$$
    \ceil{\frac{E_u X^{t+1} (\mathbf{1}\otimes I_d)}{2}}_{k} = \ceil{\frac{X^{t+1}_{ii}+X^{t+1}_{jj}}{2}},
$$
where $(i,j)\in \mathcal{E}_k$. Therefore, it holds that:
$$
    \mathbb{E}^t\left[\ceil{\frac{E_u X^{t+1} (\mathbf{1}\otimes I_d)}{2}}_{k}\right] = \mathbb{E}^t \left[\ceil{\frac{X^{t+1}_{ii}+X^{t+1}_{jj}}{2}}\right]\succ 0,
$$
which shows that $\mathbb{E}^t[Y^{t+1}]\succ 0$.\QEDB

\bibliography{references.bib}
\bibliographystyle{IEEEtran}

\begin{IEEEbiography}
[[{\includegraphics[width=1in,height=1.25in,clip,keepaspectratio]{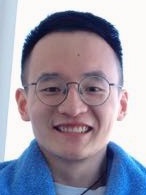}}]{Yichuan Li}
received the B.S. degree in 2016, the M.S. degree in Mechanical Engineering in 2018, the M.S. degree in Applied Mathematics, and the Ph.D. degree in Mechanical Engineering in 2022, all from the University of Illinois at Urbana-Champaign, Champaign, IL, USA. His research interests include multi-agent optimization, distributed machine learning, and control. 
\end{IEEEbiography}

\begin{IEEEbiography}
[{\includegraphics[width=1in,height=1.25in,clip,keepaspectratio]{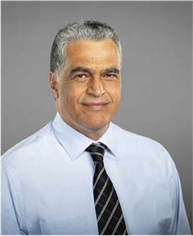}}]{Professor Petros G. Voulgaris} received the Diploma in Mechanical Engineering from the Na- tional Technical University, Athens, Greece, in 1986, and the S.M. and Ph.D. degrees in Aeronautics and Astronautics from the Massachusetts Institute of Technology in 1988 and 1991, respectively. He is currently Chair, Founding Aerospace Program Director, and Victor LaMar Lockhart Professor in Mechanical Engineering at University of Nevada, Reno. Before joining UNR in 2020 and since 1991, he has been a faculty with the Department of Aerospace Engineering, University of Illinois at Urbana-Champaign holding also appointments with the Coordinated Science Laboratory, and the depart- ment of Electrical and Computer Engineering. His research interests are in the general area of robust and optimal control and coordination of autonomous systems. Dr. Voulgaris is a recipient of several awards including the NSF Research Initiation Award, the ONR Young Investigator Award and the UIUC Xerox Award for research. He has also been a Visiting ADGAS Chair Professor, Mechanical Engineering, Petroleum Institute, Abu Dhabi, UAE and a Visiting Gaungbiao Chair at Zhejiang University, China. His research has been supported by several agencies including NSF, ONR, AFOSR, NASA. He is also a Fellow of IEEE. 
\end{IEEEbiography}

\begin{IEEEbiography}
[{\includegraphics[width=1in,height=1.25in,clip,keepaspectratio]{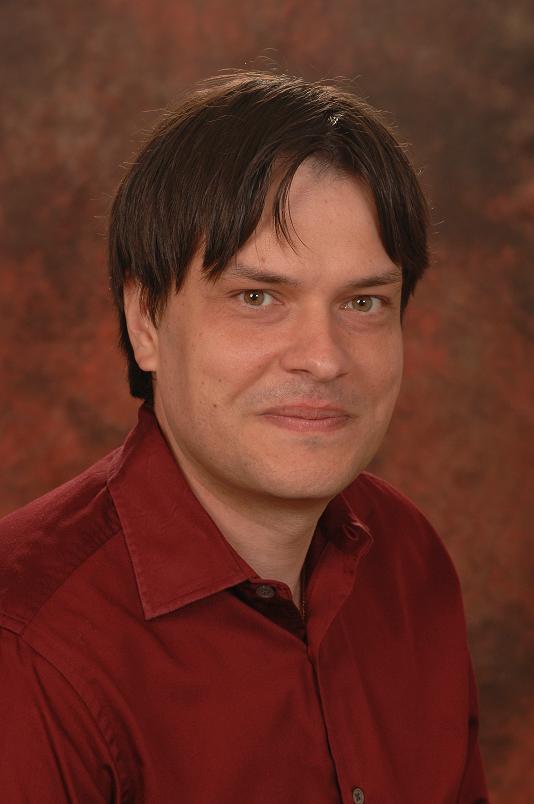}}]{Dr. Du\v{s}an M. Stipanovi\'{c}} received his B.S. degree in electrical engineering from the University of Belgrade, Belgrade, Serbia, in 1994, and the M.S.E.E. and Ph.D. degrees (under supervision of Professor Dragoslav Šiljak) in electrical engineering from Santa Clara University, Santa Clara, California, in 1996 and 2000, respectively. Dr. Stipanović had been an Adjunct Lecturer and Research Associate with the Department of Electrical Engineering at Santa Clara University (1998-2001), and a Research Associate in Professor Claire Tomlin’s Hybrid Systems Laboratory of the Department of Aeronautics and Astronautics at Stanford University (2001-2004). In 2004 he joined the University of Illinois at Urbana-Champaign where he is now Professor in the Controls Group of the Coordinated Science Laboratory and Department of Industrial and Enterprise Systems Engineering. Dr. Stipanović served as an Associate Editor on the Editorial Boards of the IEEE Transactions on Circuits and Systems I and II. Currently he is an Associate Editor for Journal of Optimization Theory and Applications.
\end{IEEEbiography}

\begin{IEEEbiography}
[{\includegraphics[width=1in,height=1.25in,clip,keepaspectratio]{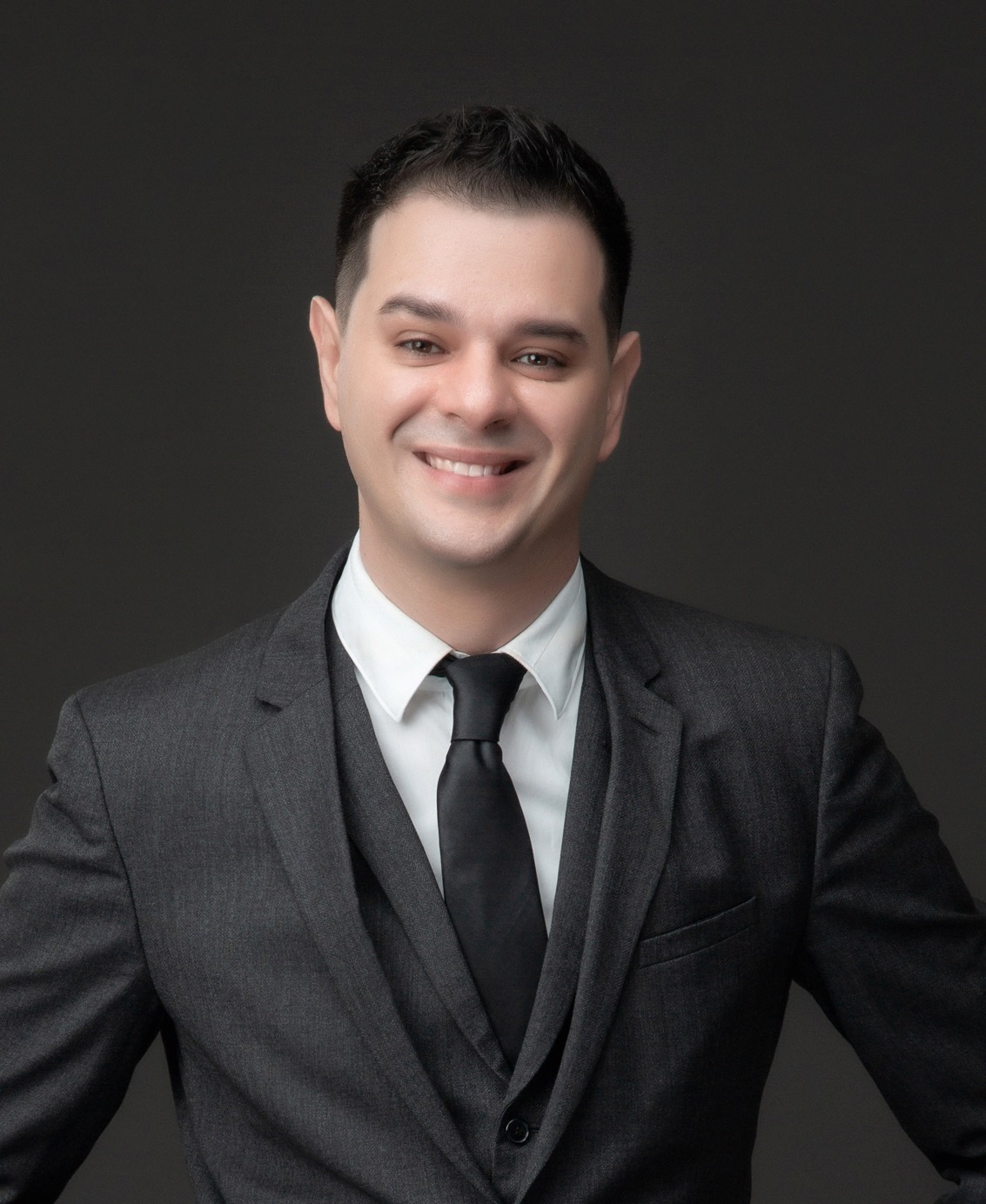}}]{Nikolaos M. Freris (Senior Member, IEEE)}
 received the Diploma in ECE from the National Technical University of Athens (NTUA), Athens, Greece, in 2005, the M.S. degree in ECE, the M.S. degree in Mathematics, and the Ph.D. degree in ECE, all from the University of Illinois at Urbana-Champaign (UIUC), Champaign, IL, USA, in 2007, 2008, and 2010, respectively. He is a Professor with the School of Computer Science and Technology and the Vice Dean of the International College at the University of Science and Technology of China (USTC), Hefei, China. His research lies in AIoT/CPS/IoT: machine learning, distributed optimization, data mining, wireless networks, control, and signal processing, with applications in power systems, sensor networks, transportation, cyber security, and robotics. Dr. Freris has published several papers in high-profile conferences and journals held by IEEE, ACM, and SIAM, and he holds three patents. His research has been sponsored by the Ministry of Science and Technology of China, Anhui Dept. of Science and Technology, Tencent, and NSF, and was recognized with the USTC Alumni Foundation Innovation Scholar award, the IBM High Value Patent award, two IBM invention achievement awards, and the Gerondelis foundation award. Previously, he was with the faculty of NYU and, before that, he held senior researcher and postdoctoral researcher positions at EPFL and IBM Research, respectively.  Dr. Freris is a Senior Member of ACM and IEEE, and a member of CCF and SIAM.
\end{IEEEbiography}
\vfill

\end{document}